\documentclass[11pt, reqno ]{amsart}
\usepackage{geometry}                
\usepackage{graphicx}
\usepackage{amssymb}
\usepackage{epstopdf}
\usepackage{amsmath}
\usepackage{mathrsfs}
\usepackage[all]{xy}
\usepackage{color}
\usepackage{tikz}
\usepackage[colorlinks=true, pdfstartview=FitV,linkcolor=blue,citecolor=blue,urlcolor=blue]{hyperref}

    \advance \topmargin by -\headheight
    \advance \topmargin by -\headsep

    \evensidemargin \oddsidemargin
\numberwithin{equation}{section}

\begin{document}

 \newtheorem{theorem}{Theorem}[section]
  \newtheorem{prop}[theorem]{Proposition}
  \newtheorem{cor}[theorem]{Corollary}
  \newtheorem{lemma}[theorem]{Lemma}
  \newtheorem{defn}[theorem]{Definition}
  \newtheorem{ex}[theorem]{Example}
    \newtheorem{conj}[theorem]{Conjecture}

 \newcommand{\cx}{{\bf C}}
\newcommand{\la}{\langle}
\newcommand{\ra}{\rangle}
\newcommand{\res}{{\rm Res}}
\newcommand{\expp}{{\rm exp}}
\newcommand{\Sp}{{\rm Sp}}
\newcommand{\lgp}{\widehat{\rm G}( { F} ((t)) )}
\newcommand{\svee}{\scriptsize \vee}
\newcommand{\deff}{\stackrel{\rm def}{=}}
\newcommand{\cal}{\mathcal}

\title{Chiral de Rham complex on the upper half plane and modular forms  }

\author[Xuanzhong Dai]{Xuanzhong Dai}
\address{Shanghai Center For Mathematical Sciences, Fudan University, Jiangwan Campus, 2005 Songhu Road, Shanghai, 200438, China}
\email{xzdai@fudan.edu.cn}

\thanks{}

\maketitle

\maketitle

\begin{abstract}
For any congruence subgroup $\Gamma$, we study the vertex operator algebra $\Omega^{ch}(\mathbb H,\Gamma)$ constructed from the $\Gamma$-invariant global sections of the chiral de Rham complex on the upper half plane, which are holomorphic at all the cusps.
We introduce an $SL(2,\mathbb R)$-invariant filtration on the global sections and show that the $\Gamma$-invariants on the graded algebra is isomorphic to certain copies of modular forms.
We also give an explicit formula for the lifting of modular forms to $\Omega^{ch}(\mathbb H,\Gamma)$ and compute the character formula of $\Omega^{ch}(\mathbb H,\Gamma)$.
Furthermore, we show that the vertex algebra structure modifies the Rankin-Cohen bracket, and the modified bracket becomes non-zero between constant modular forms involving the Eisenstein series.
\end{abstract}

\section { Introduction }

The chiral de Rham complex $\Omega_X^{ch}$ constructed by Malikov, Schectman and Vaintrob \cite{MSV} in 1998,
 is a sheaf of vertex algebras on any nonsingular algebraic variety or complex manifold $X$. 
For any open subset $U$ of $X$, $\Omega^{ch}_X(U)$ has
 a topological vertex algebra structure given by four distinguished elements, namely the Virasoro element $\omega$, an even element $J$,
 and two odd elements $Q$ and $G$ (see in \cite{MSV}).  
  It is endowed with a $\mathbb Z$-gradation from the semisimple operator $J_0$, whose eigenvalues are called the fermionic charge.
 And $\Omega^{ch}_X(U)$ together with the chiral differential $d=-Q_{0}$ becomes a complex which is quasi-isomorphic to its conformal weight zero part, namely the usual de Rham complex.

We will study the global sections of $\Omega^{ch}_X$. When $X$ is a projective $n$-space, it is a module over $\widehat{\mathfrak{sl}}_{n+1}$ \cite{MS}. When $X$ is a $K3$ surface, it is isomorphic to the simple $\mathcal N=4$ algebra with central charge $c=6$ \cite{S1}. When $X$ is a compact Ricci-flat K\"{a}hler manifold, the global sections can be viewed as an invariant subspace under the action of certain Lie algebra \cite{S2}.
In this paper, we focus on the chiral de Rham complex on the upper half plane $\mathbb H$, and denote the global sections by
 \[
 \Omega^{ch}(\mathbb H)=(V_1 \otimes \bigwedge_1 \nolimits ) \otimes_{\mathbb C[b_0]} \mathcal O(\mathbb H),
 \]
 where $V_1$ is the vacuum representation of the Heisenberg Lie algebra with the basis $a_{n},b_{n} (n\in \mathbb Z)$ and the center $C_1$ with the relation that
 \[
 [a_n,b_m]=\delta_{n,-m}C_1,
 \]
 and $\bigwedge_1$ is the vacuum representation of the Clifford algebra with the basis $\psi_n,\phi_n (n\in \mathbb Z)$ and the center $C_2$ with the relation that
 \[
 [\phi_n,\psi_m]=\delta_{n,-m}C_2.
 \]
  As vector spaces, $V_1$ is a symmetric algebra generated by $a_{-n},b_{-m}$ for $n\geq 1,m\geq 0$, and $\bigwedge_1$ is an exterior algebra generated by $\psi_{-n},\phi_{-m}$ for $n\geq 1,m\geq 0$.
Let $\Gamma \subset SL(2,\mathbb R)$ be an arbitrary congruence subgroup. 
We consider the $SL(2,\mathbb R)$-action on $\Omega^{ch}(\mathbb H)$ induced from the fractional linear transformation on $\mathbb H$ as automorphisms of vertex algebras, 
and consider the subspace $\Omega^{ch}(\mathbb H,\Gamma)$ consisting of $\Gamma$-invariant elements that are holomorphic at the cusps (see Section \ref{Section 3.3} for a detailed definition).
 Since the $SL(2,\mathbb R)$-action preserves the conformal weights of $\Omega^{ch}(\mathbb H)$, $\Omega^{ch}(\mathbb H,\Gamma)$ is naturally a $\mathbb Z_{\geq 0}$-graded vertex operator algebra.

We introduce a decreasing filtration $\{W_n\}$ of free $\mathcal O(\mathbb H)$-modules on $\Omega^{ch}(\mathbb H)$,
 which is preserved under the action of $SL(2,\mathbb R)$.
The graded algebra $Gr \Omega^{ch}(\mathbb H)=\oplus_{n\in \mathbb Z} W_n/W_{n+1}$ is still a vertex operator algebra, equipped with a much simpler induced infinitesimal action of $\mathfrak{sl}_2$.
 For any integer $k$ and $l$, we consider the invariant subspace of $W_n$, spanned by elements of conformal weight $k$ and fermionic charge $l$, denoted by $W_n(k,l)$. 
Let $W_n(k,l)^\Gamma_0$ be the subspace of $W_n(k,l)^\Gamma$ spanned by elements that are holomorphic at the cusps (see in Section \ref{Section 3.3}), and denote by $M_n(\Gamma)$ the space of modular forms of weight $n$ with respect to $\Gamma$.
Then the successive quotient $W_n(k,l)^\Gamma_0/W_{n+1}(k,l)^\Gamma_0$ can be easily embedded into $N$ copies of $M_n(\Gamma)$, where $N$ is the dimension of $ W_n(k,l)/ W_{n+1}(k,l)$ as an $\mathcal O(\mathbb H)$-module.
Our first main result in this paper is that (Theorem \ref{theorem3.5})
 \begin{equation} \label{1.1}
 W_n(k,l)^\Gamma_0/ W_{n+1}(k,l)^\Gamma_0 \cong M_{2n}(\Gamma)^{\oplus N }.
 \end{equation} 
 Note that the parallel result of (\ref{1.1}) about the chiral differential operators is proved in \cite{D}, and moreover our improved method gives explicit formulas for invariant global sections.

For any partition $\lambda=(\lambda_1,\cdots,\lambda_d)$, we define $p(\lambda):=d$, and for any symbol $X=a,b,\psi$, we denote by  $X_{-\lambda}$ the expression $X_{-\lambda_1}\cdots X_{-\lambda_d}$, 
and let $\phi_{-\lambda}:=\phi_{-\lambda_1+1}\cdots \phi_{-\lambda_{d}+1}$.
 We will show that

\begin{theorem}\label{theorem1.1} 
Let $(\lambda_0,\mu_0,\nu_0,\chi_0)$ be a four-tuple of partitions with $-p(\lambda_0)+p(\mu_0)-p(\nu_0)+p(\chi_0)=n_0$ and $f$ be a modular form of weight $2n_0$ with respect to $\Gamma$.
\begin{enumerate}
\item When $n_0>0$, then 
\begin{equation}\label{1.2}
\sum_{n=0}^\infty \frac{(2n_0-1)!}{n!(n+2n_0-1)!} D^{n}(a_{-\lambda_0}\phi_{-\mu_0} \psi_{-\nu_0} b_{-\chi_0}) f^{(n)} (b)\in \Omega^{ch}(\mathbb H,\Gamma),
\end{equation}
 is a lifting of $f$ in $\Omega^{ch}(\mathbb H,\Gamma)$.
 \item When $n_0=0$, then $f$ is a constant function, and
 \begin{equation}\label{1.3}
a_{-\lambda_0}\phi_{-\mu_0} \psi_{-\nu_0} b_{-\chi_0} +\frac{\pi i}{6}\sum_{n=1}^\infty \frac{1}{n!(n-1)!}  D^{n}(a_{-\lambda_0}\phi_{-\mu_0} \psi_{-\nu_0} b_{-\chi_0}) E_2^{(n-1)}(b)\in\Omega^{ch}(\mathbb H,SL(2,\mathbb Z)),
\end{equation}
is a lifting of $1$ in $\Omega^{ch}(\mathbb H,SL(2,\mathbb Z))$.
 \end{enumerate}
\end{theorem}
\noindent \emph{Remark.} The summations in (\ref{1.2}) and (\ref{1.3}) are finite, as the operator $D$ introduced in Section \ref{Section 4.1} is nilpotent. We call (\ref{1.2}) when $n_0>0$ (resp.  (\ref{1.3}) multiplied by the scaler function $f$ when $n_0=0$) the lifting formula of $f$ with the leading term $a_{-\lambda_0}\phi_{-\mu_0} \psi_{-\nu_0} b_{-\chi_0} f(b)\in W_{n_0}$. 
Therefore the lifting formula for a nonconstant modular form $f$ is obtained by applying invariant vertex operators to $f$, 
and the lifting of a constant modular form is obtained by invariant vertex operators and the quasi-modular form $E_2$.
Such lifting is essentially unique, modulo liftings of modular forms of higher weight.

We have the following character formula.

\begin{theorem} \label{theorem5.3.3}
The character formula of $\Omega^{ch}(\mathbb H,\Gamma)$ is given by
\begin{equation} \label{1.4}
\sum_{m,n=0}^\infty \sum_{u=0}^n \sum_{v=0}^{m+n}\dim M_{2m}(\Gamma)q^{m+2n+\frac{1}{2}u(u-1)+\frac{1}{2} v(v-3)}  \prod_{i=1}^u \dfrac{1}{1-q^i} \prod_{j=1}^v \dfrac{1}{1-q^j} \prod_{k=1}^{n-u} \dfrac{1}{1-q^k}  \prod_{l=1}^ {m+n-v}\dfrac{1}{1-q^l}.
\end{equation}
\end{theorem}

 We also show that $\Omega^{ch}(\mathbb H,\Gamma)$ is again a topological vertex algebra. The Virasoro element $\omega$ and the element $G$ are invariant under $SL(2,\mathbb R)$-action and hence they are still contained in $\Omega^{ch}(\mathbb H,\Gamma)$. But the element $J$ and $Q$ are not fixed by $\Gamma$ in general.
So we replace $J$ and $Q$ by the lifting of the constant modular form with the leading term $J$ and $Q$ respectively, namely
\begin{gather}\label{1.5}
\tilde J=J+\frac{\pi i}{3}b_{-1 }E_2(b)\\
\tilde Q=Q -\frac{\pi i}{3}\phi_{-1} E_2(b)-\frac{\pi i}{3}\phi b_{-1} E'_2(b).\label{1.6}
\end{gather}
Then the fields corresponding to $\omega$, $\tilde J$, $\tilde Q$ and $G$
make $\Omega^{ch}(\mathbb H,\Gamma) $ a topological vertex algebra. And $\Omega^{ch}(\mathbb H,\Gamma)$ equipped with the chiral differential $\tilde Q_0=Q_0$ is again a complex.

The later part of this paper is to discuss the relations between the invariant global sections of the chiral de Rham complex on the upper half plane and the Rankin-Cohen bracket. 
The Rankin-Cohen bracket defined by Cohen in 1977,  is a family of bilinear operations which send two modular forms to another modular form. 
For any modular forms $f_1$ of weight $k$, and $f_2$ of weight $l$, the $n$-th$ (n\geq 0)$ Rankin-Cohen bracket of $f_1$ and $f_2$, denoted by $[f_1,f_2]_n$ is a modular form of weight $k+l+2n$.
And the properties of the Rankin-Cohen bracket and its relations with Jacobi forms and pseudodifferential operators are studied intensively in \cite{Z2} and \cite{CMZ}. 
However the Rankin-Cohen bracket $[1,1]_n$ for constant modular forms are all $0$ for $n\geq 1$ and $[1,1]_0=1$.
In this paper, we give a natural modification of the Rankin-Cohen bracket thanks to the vertex algebra structure of $\Omega^{ch}(\mathbb H,\Gamma)$, and the formulas are non-zero with
 the quasi-modular form $E_2(b)$ involved. And the modified bracket
\[
[1,1]^\thicksim_n:= \frac{1}{144 (2\pi i)^{n-2}}\sum_{\substack{r+s=n-2}} (-1)^{r+1} {n-1\choose r}{n-1 \choose s} E_2^{(r)}(\tau) E_2^{(s)}(\tau)+\dfrac{(-1)^n+1}{12n (2\pi i)^{n-1 }} E_2^{(n-1)} (\tau)
\]
is a modular form of weight $2n$ with respect to $SL(2,\mathbb Z)$.
By the isomorphism (\ref{1.1}), we can choose a linear basis consisting of the liftings of the form (\ref{1.2}) and (\ref{1.3}),
and use the modified Rankin-Cohen bracket to describe the vertex operators in $\Omega^{ch}(\mathbb H,\Gamma)$.

 Our explicit formulas (\ref{1.2}) and (\ref{1.3}) allow us to show that $\Omega^{ch}(\mathbb H, \Gamma)$ is filtered by a chain of vertex algebra ideals

\[
L_0= \Omega^{ch}(\mathbb H,\Gamma) \supset L_1 \supset L_2 \supset \cdots
\]
 with the simple maximum quotient $\Omega^{ch}(\mathbb H,\Gamma)/L_1$,
 where $L_i$ is spanned by the liftings of modular forms of weight greater than or equal to $i$.

In the last part of the paper, we discuss the Hecke operator action on $\Omega^{ch}(\mathbb H,\Gamma)$.
We extend the group action on $\Omega^{ch}(\mathbb H)$ to the group $GL(2,\mathbb R)_{>0}$, and introduce the Hecke operator on  $\Omega^{ch}(\mathbb H,\Gamma)$.
We conclude that the Hecke algebra action and the lifting are commutative (Proposition \ref{prop7.1}).

\

The structure of this thesis is as follows. 
In Section \ref{Section 2}, we review the definition of vertex algebras and the construction of the chiral de Rham complex in \cite{MSV}. In Section \ref{Section 3}, we introduce an $SL(2,\mathbb R)$-invariant filtration on the global section of the chiral de Rham complex on the upper half plane, 
and explain the relations between $\Gamma$-fixed points on the graded vertex algebra associated to the filtration and modular forms. In Section \ref{Section 4}, we prove the isomorphism (\ref{1.1}) by showing (\ref{1.2}) and (\ref{1.3}). In Section \ref{Section 5} we calculate the cohomology groups and the character formula of $\Omega^{ch}(\mathbb H,\Gamma)$.
 In Section \ref{Section 6}, we use the generalized Cohen-Kuznetsov lifting of the constant modular form to modify the Rankin-Cohen bracket, 
 and study the structure of the vertex algebra $\Omega^{ch}(\mathbb H,\Gamma)$.

\section{Recollections}\label{Section 2}

\subsection{Vertex algebras}
\begin{defn}
A vertex operator (super)algebra is a superspace $V$, equipped with a vector $1\in V_{\bar{0}}$, a parity preserving linear map (called the state-field correspondence) from $V$ to $ End\,V[[z,z^{-1}]]$,
\begin{align*}
V \longrightarrow & End\,V[[z,z^{-1}]]\\
u \longmapsto & Y(u,z)=\sum_{n\in\mathbb Z} u_{(n)}z^{-n-1}
\end{align*}
a linear map $T\in( End\,V)_{\bar{0}}$, satisfying the following axioms
\begin{enumerate}
\item (the truncation condition): For every two vectors $u,v\in V$,
\begin{equation}
u_{(n)}v=0
\end{equation}
for $n$ sufficiently large.
\item (vacuum): $T1=0, Y(1,z)=id$, $Y(u,z)1\in V[[z]]$ and $Y(u,z)1|_{z=0}=u$. 
\item (translation covariance): 
\begin{equation}\label{2.2}
[T, Y(u,z)]=\partial_z Y(u,z),
\end{equation}
\item (locality): For every $u,v\in V$, $Y(u,z)$ and $Y(v,z)$ are mutually local, namely, there exists $N\in \mathbb Z_{>0}$, such that
\begin{equation}
[Y(u,z),Y(v,w)](z-w)^N=0.
\end{equation} 
\end{enumerate}
\end{defn}

And a vertex algebra is called a vertex operator algebra if there is a distinguished vector $\omega$ (called the Virasoro element), such that $L_0$ is semisimple and 
\[
L_{-1}=T, \; [L_m,L_n]=(m-n)L_{m+n}+\dfrac{m^3-m}{12}\delta_{m,-n}c
\]
where $L_n=\omega_{(n+1)}$, and $c\in \mathbb C$ is a constant called the central charge. 

For any $a,b\in V$, we will frequently use the Borcherds identity
\begin{equation} \label{2.4}
(a_{(n)}b)_{(m)}=\sum_{i=0}^{\infty} (-1)^i {n\choose i}a_{(n-i)}b_{(m+i)}- \sum_{i=0}^\infty (-1)^{p(a)p(b)+n+i}{n\choose i} b_{(m+n-i)} a_{(i)},
\end{equation}
where we denote by $p:V\to \{0,1\}$ the parity function.

Below we will briefly review two basic examples of vertex operator algebras.

\begin{ex}
 We fix a positive integer $N$. Let $\mathcal H$ be the Heisenberg Lie algebra with the basis $a^i_{n},b^i_{n}, i=1,\cdots,N, n\in \mathbb Z$, the central element $C$, and the nontrivial commutation relations
\begin{equation}
[a^i_m,b^j_n]=\delta_{i,j} \delta_{n+m,0}C.
\end{equation}
The Heisenberg Vertex algebra $V_N$ is defined to be the vacuum representation of $\mathcal H$, which as an induced module of $\mathcal H$, is generated by the vacuum vector $1$, with the following relations
\begin{equation}
a^i_m 1=0 \;\;\text{ if }m\geq 0; \;\;b^i_n1=0\;\;\text{ if } n\geq 1;\;\; C1=1.
\end{equation} 
As a vector space, $V_N$ is a symmetric algebra generated by elements $a^i_{m},b^i_n, m<0,n\leq 0, 1\leq i \leq N$. For $a^i=a^i_{-1}\cdot 1$, and $b^i=b^i_{0}\cdot 1$, the corresponding fields are
\begin{align*}
Y(a^i,z)&=a^i(z)=\sum_{n\in \mathbb Z} a^i_n z^{-n-1},\\
Y(b^i,z)&=b^i(z)=\sum_{n\in \mathbb Z} b^i_n z^{-n},
\end{align*}
and the nontrivial operator product expansion of the basic fields is
\begin{equation}
a^i(z)b^j(w)\sim \delta_{ij} (z-w)^{-1}.
\end{equation}

Define 
\begin{equation}
\omega:=\sum_{i=1}^N b^i_{-1}a^i_{-1}\cdot 1.
\end{equation}
Then $\omega$ is a Virasoro element with central charge $2N$.
So $V_N$ is a vertex operator algebra and it has the conformal field $L(z)=\sum_{i=1}^N :\partial_z b^i(z) a^i(z):$ with OPE
\begin{equation}
L(z)L(w)\sim \frac{N}{(z-w)^{4} }+\frac{2L(w) }{(z-w)^{2}}+\frac{\partial_w L(w) }{z-w}.
\end{equation}
\end{ex}

\begin{ex}
Let $Cl$ be the Clifford algebra with the basis consisting of odd elements $\phi^i_n, \psi^i_n, i=1,\cdots, N, n\in \mathbb Z$, and the central element $C$, and the nontrivial commutation relations
\begin{equation}
[\phi^i_{m}, \psi^j_{n}]=\delta_{i,j} \delta_{m+n,0}.
\end{equation}
The Clifford vertex algebra $\bigwedge_N$ is defined to be the vacuum representation of $Cl$, generated by the vacuum vector $1$, and the relations
\[
\phi^i_{m} 1=0  \;\;\text{ if } m\geq 1; \;\; \psi^i_n 1=0 \;\;\text{ if } n\geq 0; \;\; C1=1.
\]
As a vector space, $\bigwedge_N$ is the exterior algebra generated by elements $\phi^i_m ,\psi^i_n$ for $m\leq 0, n<0, 1\leq i\leq N$. Let $\phi^i=\phi^i_{0}\cdot1$, and $\psi=\psi^i_{-1}\cdot 1$. The basic odd fields are given as follows
\begin{align*}
Y(\phi^i,z)&=\phi^i(z) =\sum_{n\in\mathbb Z} \phi^i_n z^{-n},\\
Y(\psi^i,z)&=\psi^i(z) =\sum_{n\in \mathbb Z} \psi^i_nz^{-n-1}.
\end{align*}
Then $\bigwedge_N$ is a vertex operator algebra with the Virasoro element 
\begin{equation}
\omega:=\sum_{i=1}^N \phi^i_{-1}\psi^i_{-1}
\end{equation}
of central charge $-2N$ and $L(z)=\sum_{i=1}^N :\partial_z \phi^i(z)\psi^i(z):$ with OPE
\begin{equation}
L(z)L(w)\sim \frac{-N}{(z-w)^4}+\frac{2L(w)}{(z-w)^2}+\frac{\partial_w L(w)}{z-w}.
\end{equation}
\end{ex}

\subsection{Chiral De Rham Complex} \label{Section 2.2}
 We consider the tensor product vertex algebra
\begin{equation}
\Omega_N:=V_N\otimes \bigwedge \nolimits_N,
\end{equation}
where the Virasoro element is given by $\omega=\sum _{i=1}^N b^i_{-1}a^i_{-1}+\phi^i_{-1}\psi^i_{-1}$ with central charge $0$. There are three other special elements in $\Omega_N$, namely an even element $J:=\sum_{i=1}^N \phi^i \psi^i_{-1}$ of conformal weight $1$, two odd elements $Q=\sum_{i=1}^N a^i_{-1} \phi^i$ and $G=:\sum_{i=1}^N \psi^i_{-1} b^i_{-1}$ of conformal weight $1$ and $2$ respectively. The fields corresponding to the above four elements satisfy the following OPEs: \cite{MSV}
\begin{gather}
\label{2.14}L(z)L(w)\sim \frac{2L(w)}{(z-w)^{2}}+\frac{L(w)'}{z-w},\; L(z)J(w)\sim \frac{-N}{(z-w)^{3}} +\frac{J(w)}{(z-w)^2}+\frac{J(w)'}{z-w},ct\\
L(z)Q(w) \sim \frac{Q(w)}{(z-w)^2}+\frac{Q(w)'}{z-w},\; L(z)G(w)\sim \frac{2G(w)}{(z-w)^w} +\frac{G(w)'}{z-w} ,\\
J(z) J(w)\sim \frac{N}{(z-w)^2},\;J(z)Q(w)\sim \frac{Q(w)}{z-w}, \;J(z)G(w)\sim \frac{-G(w)}{z-w},\\
\label{2.17}Q(z)Q(w)\sim 0,\; Q(z) G(w) \sim \frac{N}{(z-w)^3} +\frac{J(w)}{(z-w)^2}+\frac{L(w)}{z-w},\; G(z)G(w)\sim 0 .
\end{gather}
The four fields $L(z), J(z),Q(z)$ and $G(z) $ with relations (\ref{2.14})-(\ref{2.17}) make $\Omega_N$ into a topological vertex algebra of rank $N$.
Let $V'_N$ be the subalgebra of $V_N$ generated by $a^i_{-1}$ and $b_{-1}^i$ for $i=1,\cdots,N$, and 
\begin{equation}
\Omega'_N= V'_N \otimes \bigwedge \nolimits _N,
\end{equation}
is again a topological vertex algebra.

We denote by $\Omega_N(m)\subset \Omega_N$ the subspace consisting of elements with fermionic charge $m$, i.e.
\begin{equation}
\Omega_N(m):=\{ w\in \Omega_N \;| \; J_{(0)} w=mw\},
\end{equation} 
where $J_{(0)}$ is called the fermionic charge operator, and it acts semisimply on $\Omega_N$. The fermionic charge equals $0$ for the generators $a^i,b^i$, and  $1,-1$ for $\phi^i, \psi^i$ respectively. 
According to the first OPE in (\ref{2.17}), we have
\begin{equation}
Q_{(0)}Q=0.
\end{equation}
Therefore the space $\Omega_N=\oplus_{m\in \mathbb Z} \Omega_N(m)$ together with the chiral de Rham differential $d:= -Q_{(0)}$, becomes a complex.

\

Let $U$ be an open subset of an $n$ dimensional complex manifold (or smooth algebraic variety), with local coordinates $b^{1},\cdots , b^{N}$. Let $\mathcal O(U)$ be the space of smooth functions (or algebraic functions ) on $U$. 
Then $\mathcal O(U)$ is a $\mathbb C[b^1_{0},\cdots,b^N_0]$-module, where the action of $b^i_0$ on $\mathcal O(U)$ is simply the multiplication by $b^i$. Then $\Omega^{ch}(U)$ is defined to be the localization of $\Omega_N$ on $U$, namely
\[
\Omega^{ch}(U)=\Omega_N \otimes _{\mathbb C[b^1_0,\cdots,b^N_0]} \mathcal O(U).
\]
Then $\Omega^{ch}(U)$ is a vertex algebra generated by $a^i(z),\partial b^i(z), \phi^i(z),\psi^i(z)$ and $Y(f,z)$ for  $f\in \mathcal O(U)$, where the field $Y(f,z)$ is defined by
\begin{equation} \label{2.21}
Y(f,z)=\sum_{i= 0}^\infty \frac{\partial^i}{i!}  f(b) (\sum_{n\neq 0} b_nz^{-n})^i.
\end{equation}
We write $f(b)_{m+1}:=f(b)_{(m)}$ for the coefficient of $z^{-m-1}$ in the field $Y(f,z)$. And these generators satisfy the following nontrivial OPEs,
\begin{gather*}
a^i(z)\partial b^j(w)\sim \frac{\delta_{ij} }{(z-w)^2},\; \psi^i(z)\phi^j(w)\sim \frac{\delta_{ij}}{z-w},\\
a^i(z)f(w)\sim \frac{\frac{\partial}{\partial \, b^i} f(w)}{z-w}.
\end{gather*}

If we have another coordinates $\tilde b^1,\cdots,\tilde b^N$ on $U$,
the coordinate transformation equations for the generators are
\begin{align*}
\tilde a^i_{-1} =& a^j_{-1} \frac{\partial b^j}{\partial \tilde b^i}  +  \frac{\partial^2 b^j}{\partial \tilde b^i \partial \tilde b^m}  \frac{\partial \tilde b^m }{\partial b^r} \phi^r \psi^j,\\
\tilde b^i_{-1}= & \frac{\partial \tilde b^i}{\partial b^j} b^j_{-1},\\
\tilde \phi^i_0\;= &\frac{\partial \tilde b^i}{\partial b^j} \phi^i _0,\\
\tilde \psi^i_{-1}=& \frac{\partial b^j}{\partial \tilde b^i} \psi^j_{-1},
\end{align*}
where we use Einstein summation convention.

\section{Invariant Global Sections} \label{Section 3}

In this section, we apply the construction of the chiral de Rham complex in \cite{MSV} to the upper half plane. 
And we will introduce an $SL(2,\mathbb R)$-action on the global sections, and consider the invariant sections under the induced action of congruence subgroups of the modular group $SL(2,\mathbb Z)$.

\subsection{Localization on the Upper Half Plane}

From now on we will focus on the upper half plane $\mathbb H$
\[
\mathbb H:=\{ \tau \in\mathbb C \; | \text{ im }\tau >0\}.
\]
 In this case $N=1$, we denote by $a=a^1,b=b^1,\phi=\phi^1, \psi=\psi^1$. 
 Under the identification of $b_0$ with the variable $\tau$, we can view the polynomial ring $\mathbb C[b_0]$ as a subring of holomorphic functions on $\mathbb H$. 
According to the construction in \cite{MSV}, the global sections
\begin{equation} \label{3.1}
\Omega^{ch}(\mathbb H) := \Omega_1 \otimes _{\mathbb C[b_0]} \mathcal O(\mathbb H),
\end{equation}
 is a vertex operator algebra. 
 $\Omega^{ch}(\mathbb H)$ is generated by the basic fields $a(z),\partial b(z),\phi(z),\psi(z)$ and $Y(f,z)$ for $f\in \mathcal O(\mathbb H)$.

Certain vertex operators on $\Omega^{ch}(\mathbb H)$ generates representations of affine Kac-Moody algebra $\widehat{\mathfrak{sl}}_2$. More precisely let
\begin{equation}\label{3.2}
E:=-a_{-1},\; F:= a_{-1}b_0^2 +2b_0 \phi_0 \psi_{-1},\; H:= -2a_{-1}b_0-2\phi_0 \psi_{-1}.
\end{equation}
We have the following theorem
\begin{theorem} \label{theorem3.1} \emph{(\cite{W}, \cite{FF})}
The coefficients of $E_{(n)}, F_{(n)}, H_{(n)}$ of fields $Y(E,z), Y(F,z), Y(H,z)$ satisfy the relations of affine Kac-Moody algebra $\widehat{\mathfrak{sl}}_2$ of level $0$, where $E,F,H$ corresponds to matrices 
\[\begin{pmatrix}  0 & 1 \\ 0 & 0\end{pmatrix},\;
 \begin{pmatrix} 0 & 0\\ 1 & 0 \end{pmatrix},\;
  \begin{pmatrix} 1 & 0 \\ 0 & -1 \end{pmatrix}
  \] respectively.
\end{theorem}

Note that $E_{(0)},F_{(0)}$ and $H_{(0)}$ give an action of Lie algebra $\mathfrak{sl}_2$ on $\Omega^{ch}(\mathbb H)$ as derivations and they can be integrated to an $SL(2,\mathbb R)$-action as automorphisms of vertex algebra \cite{MSV}.
Because we will consider the action of a congruence subgroup $\Gamma \subset SL(2,\mathbb R)$, and it will be related to the theory of modular forms, where the group action is always from the right, 
we will make our action of $SL(2,\mathbb R)$ a right action. 
By definition, for $g=e^{x}\in SL(2,\mathbb R),\; x\in \mathfrak{sl}_2$, then
\[
\pi(g)=\sum_{n\geq 0} \frac{(-x_{(0)})^n}{n!}.
\]
And we have $\pi(g_1g_2)=\pi(g_2)\pi(g_1)$.
The $SL(2,\mathbb R)$-action commutes with the translation operator $T=L_{-1}=\omega_{(0)}$ for the fact that
\[
[T,x_{(0)}]=(Tx)_{(0)}=0,\;\;\;\text{ for } x\in\mathfrak{sl}_2\subset \Omega^{ch}(\mathbb H).
\]
 And it also commutes with the semisimple operator $L_0=\omega_{(1)}$, so it preserves the conformal weight.

The formulas of the action of 
 \begin{equation} \label{3.3}
 g=\begin{pmatrix}
 \alpha  & \beta\\
 \gamma &\delta
 \end{pmatrix}\in SL(2,\mathbb R)
 \end{equation} 
on the generators $a,b_{-1},\psi,\phi,f(b)\in \Omega^{ch}(\mathbb H)$ are given as follows
\begin{align} 
\nonumber \pi(g) a&=a_{-1}(\gamma b+\delta)^2 +2 \gamma (\gamma b+\delta) \phi _0 \psi_{-1},\\
\nonumber\pi(g) b_{-1}&=b_{-1} (\gamma b+\delta)^{-2},\\
\label{3.4}\pi(g) \psi &=\psi_{-1} (\gamma b+\delta)^{2},\\
\nonumber\pi(g) \phi &=\phi_0 (\gamma b+\delta)^{-2},\\
\nonumber\pi(g)f(b)&=f(gb)=f\left(\dfrac{\alpha b+\beta}{\gamma b +\delta} \right).
\end{align}
Note that (\ref{3.4}) agrees with the coordinate transformation equations in Secction \ref{Section 2.2}.

For $x\in \Omega^{ch}(\mathbb H)$, the adjoint action of $g\in SL(2,\mathbb R)$ on the operator $x_{(n)}$ is defined to be
\begin{equation}\label{3.5}
\pi(g) x_{(n)} \pi(g)^{-1}=(\pi(g)x)_{(n)}.
\end{equation}
We compute the adjoint action on $a_n$ as an example. We take $x=a$ in (\ref{3.5}), and apply the action formula of $\pi(g)a$ in (\ref{3.4}). Then 
\[
\pi(g)a_{n} \pi(g)^{-1}=(\pi(g)a)_{(n)} =(a_{-1} (\gamma b+\delta)^2)_{(n)} +2\gamma((\gamma b+\delta) \phi \psi)_{(n)}.
\]
Using Borcherds identity (\ref{2.4}), the right hand side equals
\[  \sum_{k\geq 1} a_{-k} (\gamma b+\delta)^2 _{n+k} + \sum_{k\geq 0} (\gamma b+\delta)^2_{n-k}a_k  +2 \gamma\sum_k  (\gamma b+\delta)_{-k}\left(\sum_{i\geq 0} \phi_{-i} \psi_{n+k+i}  -\sum_{i\geq 1} \psi_{n+k-i}\phi_i  \right) ,
\]
where $(\gamma b+\delta)^i_m=(\gamma b+\delta)^i_{(m-1)}$ is the $(m-1)$th Fourier coefficients of the field $Y((\gamma b+\delta)^i, z)$ as in (\ref{2.21}). 
Notice that (\ref{2.21}) is equivalent to the following identities
\begin{equation}\label{3.6}
f(b)_{0}=f(b)_{(-1)}= f(b)+ \sum_{l=2}^\infty \frac{\partial^l}{l!} f(b) \sum_{\substack{i_1,\cdots,i_l   \neq 0 \\ i_1+\cdots+i_l=0    }  }  b_{i_1}\cdots b_{i_l},
\end{equation}
 and for $k\neq 0$,
 \begin{equation}\label{3.7}
 f(b)_{k}=f(b)_{(k-1)}=\sum_{l\geq 1} \frac{\partial^l}{l!} f(b) \sum_{\substack{i_1,\cdots,i_l \neq 0 \\ i_1+\cdots+i_l=k  }} b_{i_1}\cdots b_{i_l}.
 \end{equation}
As a consequence, we have
\begin{align}
\pi(g) a_n \pi(g)^{-1}\nonumber
&=\sum_{k\geq 1} a_{-k } \left(\delta_{k,-n} (\gamma b+\delta)^2 +(1-\delta_{k,-n})2\gamma (\gamma b+\delta) b_{n+k} +\gamma ^2 \sum_{\substack{i,j\neq 0\\ i+j=n+k}} b_{i}b_j       \right)\\
\nonumber&+\sum_{k\geq 0}\left(\delta_{k,n} (\gamma b+\delta)^2 +(1-\delta_{k,n})2\gamma (\gamma b+\delta) b_{n-k} +\gamma ^2 \sum_{\substack{i,j\neq 0\\ i+j=n-k}} b_{i}b_j       \right) a_k\\
\label{3.8}&+2\gamma (\gamma b+\delta)  \left(\sum_{i\geq 0}  \phi_{-i}\psi_{n+i}-\sum_{i\geq 1} \psi_{n-i} \phi_{i}        \right)               +   2\gamma^2 \sum_{k\neq0}  b_{-k} \left(\sum_{i\geq 0} \phi_{-i} \psi_{n+k+i}  -\sum_{i\geq 1} \psi_{n+k-i}\phi_i  \right),                                
\end{align}
where $\delta_{m,n}$ is the Chronecker symbol.
Similarly the formulas of the adjoint action on the operators $\phi_{n},\psi_{n}$ and $b_{n}$ for $n\in \mathbb Z$ are given by
\begin{align}
 \nonumber
 \pi(g) \phi_{n} \pi(g)^{-1}&=(\pi(g)\phi)_{(n-1)}= (\phi(\gamma b+\delta)^{-2})_{(n-1)}=\sum_{k} \phi_{-k} (\gamma b+\delta)^{-2}_{n+k}\\
\label{3.9} &=\sum_{k}\phi_{-k} \left( \delta_{k,-n}  (\gamma b+\delta)^{-2}+ \sum_{l\geq 1} \sum_{  \substack{ i_1,\cdots, i_l\in \mathbb Z_{\neq 0} \\ i_1+\cdots +i_l=n+k } } (-\gamma)^{l} (l+1) b_{i_1}\cdots b_{i_l}    \right)
\end{align}
\begin{align}
 \nonumber \pi(g) \psi_n \pi(g)^{-1} & =(\pi(g)\psi)_{(n)}=(\psi_{-1}(\gamma b+\delta)^{2})_{(n)}=\sum_{k} \psi_{-k}{ (\gamma b+\delta)^2}_{n+k}\\
  \label{3.10} &= \sum_{k} \psi_{-k} \left(\delta_{k,-n} (\gamma b+\delta)^2 +(1-\delta_{k,-n})2\gamma (\gamma b+\delta) b_{n+k} +\gamma ^2 \sum_{\substack{i,j\neq 0\\ i+j=n+k}} b_{i}b_j       \right)
 \end{align}
 \begin{align}
 \nonumber \pi(g)b_{n} \pi(g)^{-1}&=  (\pi(g)b)_{n}=\left(\dfrac{\alpha b+\beta}{\gamma b+\delta}  \right)_{n}\\
 \label{3.11}
&= \delta_{n,0} \dfrac{\alpha b+\beta}{\gamma b+\delta}+  \sum_{l\geq 1} \sum_{\substack{i_1,\cdots,i_l   \in \mathbb Z_{\neq 0}  :\\ i_1+\cdots+i_l=n }} (-\gamma)^{l-1}(\gamma b+\delta)^{-l-1} b_{i_1}\cdots b_{i_l}
 \end{align}
where as an analogy to the computation of $\pi(g)a_n\pi(g)^{-1}$, the third equalities in (\ref{3.9})-(\ref{3.10}) are given by (\ref{2.4}), and the last equalities in (\ref{3.9})-(\ref{3.11}) are given by (\ref{3.6}) and (\ref{3.7}).
Since $SL(2,\mathbb R)$ acts on $\Omega^{ch}(\mathbb H)$ as automorphisms, the action formula of $g$ on $a_{-\lambda} \phi_{-\mu}\psi_{-\nu} b_{-\chi}f(b)$ is given as follows,
\begin{equation} \label{3.12}
\pi(g) a_{-\lambda} \phi_{-\mu}\psi_{-\nu} b_{-\chi}f(b) =\prod_{i=1}^{p(\lambda)}(\pi(g) a)_{-\lambda_{(i)}}\cdot \prod_{j=1}^{p(\mu)} (\pi(g)\phi)_{-\mu_{(j)}}  \cdot \prod_{k=1}^{p(\nu)}(\pi(g))_{-\nu_{(k)}}\cdot \prod_{l=1}^{p(\chi)} (\pi(g)b)_{-\chi_{(l)}} \cdot f(gb).
\end{equation}

\subsection{Invariant Global Sections and Modular Forms}  \label{Section 3.3}

For any partition $\lambda=(\lambda_{1},\lambda_{2},\cdots, \lambda_{d})$ with $ \lambda_{1}\geq \lambda_{2}\geq \cdots \geq \lambda_{d}\geq 1$, we define  $|\lambda |:=\sum_{i=1}^d \lambda_{i}$. 
 We call $\mu=(\mu_{1},\mu_{2},\cdots,\mu_{t})$ with $\mu_{1} >\mu_{2}>\cdots>\mu_{t} \geq 1$, a partition with distinct parts.
And for convenience we call $(\lambda,\mu,\nu,\chi)$ a four-tuple of partitions, if $\lambda,\chi$ are partitions, and $\mu,\nu$ are partitions with distinct parts.
Every element in $\Omega^{ch}(\mathbb H)$ can be written as a sum of elements of type $a_{-\lambda}\phi_{-\mu} \psi_{-\nu}b_{-\chi}f(b)$ with $f(b)\in \mathcal O(\mathbb H)$. Notice that $a_{-\lambda}\phi_{-\mu}\psi_{-\nu}b_{-\chi} f(b)$ has conformal weight $|\lambda|+|\mu|+|\nu|+|\chi|-p(\mu)$.
We also consider the empty set as a partition, and set $a_{-\emptyset} =b_{-\emptyset}=\phi_{-\emptyset}=\psi_{-\emptyset}=1$, and $|\emptyset|=p(\emptyset)=0$. For any given monomial $v=a_{-\lambda}\phi_{-\mu}\psi_{-\nu} b_{-\chi} f \in \Omega^{ch}(\mathbb H)$, 
we call $-p(\lambda)+p(\mu)-p(\nu)+p(\chi)$ the part of $v$
and the corresponding four-tuple $(\lambda,\mu,\nu,\chi)$.

Let $\Gamma\subset SL(2,\mathbb Z)=\Gamma(1)$ be an arbitrary congruence subgroup.
As a subgroup of $SL(2,\mathbb R)$, $\Gamma$ also acts on $\Omega^{ch}(\mathbb H)$. We denote by $\Omega^{ch}(\mathbb H)^\Gamma$ the $\Gamma$-fixed points of $\Omega^{ch}(\mathbb H)$.
$\Omega^{ch}(\mathbb H)$ is not an interesting object as it is too big, so we consider the elements in $\Omega^{ch}(\mathbb H)^\Gamma$ satisfying the cuspidal conditions similar to the definition of modular forms of $\Gamma$.

We consider $\Gamma=\Gamma(1)$ first. Recall that $T= \begin{pmatrix} 1 & 1\\0 &1\end{pmatrix}$.
By (\ref{3.4}) $\pi(T)$ preserves the generators $a,b_{-1},\phi,\psi$, i.e.
\[
\pi(T)a=a,\; \pi(T) b_{-1}=b_{-1}, \; \pi(T) \phi=\phi,\; \pi(T) \psi=\psi,
\]
so $\pi(T)$ also preserves $a_{-n},b_{-n},\psi_{-n},\phi_{-m}$ for $n\geq 1$ and $m\geq 0$. 
And since it acts as an automorphism on $\Omega^{ch}(\mathbb H)$, we have
\[
\pi (T) \sum a_{-\lambda}\phi_{-\mu} \psi_{-\nu}b_{-\chi} f_{\lambda,\mu,\nu,\chi}(b)=\sum a_{-\lambda}\phi_{-\mu} \psi_{-\nu} b_{-\chi}f_{\lambda,\mu,\nu,\chi}(b+1).
\]
Hence $f_{\lambda,\mu,\nu,\chi}(b+1)=f_{\lambda,\mu,\nu,\chi}(b)$, and $f_{\lambda,\mu,\nu,\chi}(b)$ has a $q$-expansion at the cusp $\infty$,
 \[
 f_{\lambda,\mu,\nu,\chi}(b)= \sum_{m=-\infty}^{\infty} u_{\lambda,\mu,\nu,\chi}(m) q^m , \;\;\;\text{ where } q=e^{2\pi i b}.
 \]
 We call $v=\sum_{\lambda,\mu,\nu,\chi} a_{-\lambda} \phi_{-\mu}\psi_{-\nu} b_{-\chi} f_{\lambda,\mu,\nu,\chi}$ is holomorphic at $\infty$, if 
  for arbitrary four-tuple of partitions $(\lambda,\mu,\nu,\chi)$, we have $u_{\lambda,\mu,\nu,\chi}(m)=0$ for $m<0$. Since all the cusps $\mathbb Q\cup \{\infty\}$ are $SL(2,\mathbb Z)$-equivalent, we call $v$ is holomorphic at the cusps.
 
 For a general congruence subgroup $\Gamma$, the notion of holomorphicity at the cusp $c\in \mathbb Q\cup \{\infty\}$ needs more discussions. 
 Choose $\rho\in SL(2,\mathbb Z)$ such that $\rho (c)=\infty$.  
 Then $\pi(\rho)v=\sum a_{-\lambda'}\phi_{-\mu'} \psi_{-\nu'}b_{-\chi'}  \tilde {f}_{\lambda',\mu',\nu',\chi'}$ is invariant under $\rho^{-1} \Gamma \rho$, as the group action is a right action. 
 And since $\rho^{-1} \Gamma \rho$ contains the translation matrix $\begin{pmatrix}  1 & N \\ 0 & 1 \end{pmatrix}$,  for some positive integer $N$, $\pi(\rho)v$ is fixed by $\begin{pmatrix}  1 & N \\ 0 & 1 \end{pmatrix}$(cf. \cite{Bu} p.41-42),
 which implies that $\tilde{f}_{\lambda',\mu',\nu',\chi'}(b_0+N)=\tilde{f}_{\lambda',\mu',\nu',\chi'}(b_0)$.
 Hence $\tilde f_{\lambda',\mu',\nu',\chi'}$ has a Fourier expansion $\sum \tilde{u}_{\lambda',\mu',\nu',\chi'}(m) e^{2\pi i mb/N}$. 
We say that $v$ is holomorphic at the cusp $c$ if for arbitrary four-tuple $(\lambda',\mu',\nu',\chi')$, we have $\tilde{u}_{\lambda',\mu',\nu',\chi'}(m)=0$ for $m<0$. 
 We denote by $\Omega^{ch}(\mathbb H,\Gamma)$  the $\Gamma$-invariant vectors in $\Omega^{ch}(\mathbb H)$ that are holomorphic at all the cusps. Using (\ref{2.21}), we can prove that $\Omega^{ch}(\mathbb H,\Gamma)$ is a vertex subalgebra. 
 
 \begin{prop} 
$\Omega^{ch}(\mathbb H,\Gamma) $ is a vertex operator subalgebra of $\Omega^{ch}(\mathbb H)$.
\end{prop}

\
\

Now we will introduce a partial order on the collection of four-tuples of partitions. 
We say $(\lambda,\mu,\nu,\chi)>(\lambda',\mu',\nu',\chi')$ if $ -p(\lambda)+p(\mu)-p(\nu)+p(\chi)<-p(\lambda')+p(\mu')-p(\nu')+p(\chi')$. 
And we define a family of free $\mathcal O(\mathbb H)$-submodules:
\[
W_{m}:= Span _\mathbb C \{a_{-\lambda}\phi_{-\mu} \psi_{-\nu} b_{-\chi} f(b) \in \Omega^{ch}(\mathbb H) | -p(\lambda)+p(\mu)-p(\nu)+p(\chi)\geq m\},
\]
which gives an decreasing filtration on $\Omega^{ch}(\mathbb H)$, namely
\begin{equation}
W_n \subset W_m,\;\; \text{ if } n \geq m.
\end{equation}
 The following lemma shows that $SL(2,\mathbb R)$-action preserves the filtration.
 
 \begin{lemma}\label{lemma3.3} \cite{D}
 For $g$ as in (\ref{3.3}), and holomorphic function $f$ on $\mathbb H$, 
\begin{equation}\label{3.14}
\pi(g) a_{-\lambda} \phi_{-\mu}\psi_{-\nu} b_{-\chi} f(b) =a_{-\lambda}\phi_{-\mu} \psi_{-\nu} b_{-\chi} (\gamma b+\delta)^{-2m}f(gb)
+W_{m+1},
\end{equation}
where $m=-p(\lambda)+p(\mu)-p(\nu)+p(\chi)$ is the part of the four tuple $(\lambda,\mu,\nu,\chi)$.
\end{lemma}

Since $SL(2,\mathbb R)$-action preserves the conformal weight and fermionic charge by (\ref{3.4}), $W_{m}$ can be decomposed as a direct sum of submodules 
\[
W_m=\oplus_{k\geq 0,l \in \mathbb Z } W_m(k,l),
\]
where $W_m(k,l)\subset W_m$ is spanned by elements of conformal weight $k$ and fermionic charge $l$.
And (\ref{3.14}) will be turned into
\begin{equation}\label{3.15}
 \pi(g) a_{-\lambda} \phi_{-\mu}\psi_{-\nu} b_{-\chi} f(b) =a_{-\lambda}\phi_{-\mu} \psi_{-\nu} b_{-\chi} (\gamma b+\delta)^{-2m}f(gb)
+W_{m+1}(k,l),
\end{equation}
where $m=-p(\lambda)+p(\mu)-p(\nu)+p(\chi)$, $k=|\lambda|+|\mu|+|\nu|+|\chi|-p(\mu)$, and $l=p(\mu)-p(\nu)$.

Now notice that given a four-tuple of partitions $(\lambda_0,\mu_0,\nu_0,\chi_0)$ with the part $n_0$, there are only finitely many four-tuples $(\lambda,\mu,\nu,\chi)<(\lambda_0,\mu_0,\nu_0,\chi_0)$ 
satisfying that
\begin{align}
\label{3.16}|\lambda|+|\mu|+|\nu|+|\chi|-p(\mu)&=|\lambda_0|+|\mu_0|+|\nu_0|+|\chi_0|-p(\mu_0),\\
\label{3.17} p(\mu)-p(\nu)&=p(\mu_0)-p(\nu_0).
\end{align}
 Indeed, (\ref{3.16}) implies $|\lambda|+|\nu|+|\chi|$ is bounded from above, for the reason that  $p(\mu)\leq |\mu|$. Adding the two equations (\ref{3.16}) and (\ref{3.17}), we have
 \begin{equation}
\label{3.18}|\lambda|+|\mu|+|\nu|+|\chi|-p(\nu)=|\lambda_0|+|\mu_0|+|\nu_0|+|\chi_0|-p(\nu_0),
\end{equation}
which also implies that $|\lambda|+|\mu|+|\chi|$ is bounded from above for a similar reason. Hence $|\lambda|+|\mu|+|\nu|+|\chi|$ is bounded from above, and obviously it is bounded from below, as each term is nonnegative. 
Hence there are only finitely many four-tuples satisfying both (\ref{3.16}) and (\ref{3.17}), because the partition function is always finite, and  
the linear equation $x_1+x_2+x_3+x_4=n$, for arbitrary $n\geq 0$, has only finitely many nonnegative solutions.

Since the two spaces $W_{n_0+1}(k,l)$ and $ W_{n_0}(k,l)$ are preserved under $SL(2,\mathbb R)$-action for fixed $k\geq 0$, and $l\in \mathbb Z$, so is the quotient space $W_{n_0}(k,l)/W_{n_0+1}(k,l)$ under the induced group action.
So for any congruence subgroup $\Gamma $, we have a short exact sequence
\[
0\longrightarrow W_{n_0+1}(k,l)^\Gamma \longrightarrow W_{n_0} (k,l)^\Gamma \longrightarrow (W_{n_0}(k,l)/W_{n_0+1}(k,l))^{\Gamma}.
\]

Fix any four-tuple $(\lambda_0,\mu_0,\nu_0,\chi_0)$ with the part $n_0$, and consider an arbitrary element 
\[a_{-\lambda_0} \phi_{-\mu_0}\psi_{-\nu_0} b_{-\chi_0} f(b) \in W_{n_0}(k,l)
\]
 of conformal weight $k$ and fermionic charge $l$. By (\ref{3.15}),  $a_{-\lambda_0} \phi_{-\mu_0}\psi_{-\nu_0} b_{-\chi_0} f(b)$ is fixed by $\Gamma$ modulo $W_{n_0+1}(k,l)$, if and only if 
\begin{equation} \label{3.19}
f(b)=(\gamma b+\delta)^{-2n_0} f(gb),\;\;\text{ for any } g\in \Gamma.
\end{equation}
 Note that (\ref{3.19}) is the automorphy condition for modular forms of weight $2n_0$. We denote by $W_{n_0}(k,l)^\Gamma_0$ the subspace of $W_{n_0}(k,l)^\Gamma$ consisting of elements holomorphic at all the cusps. 
 So we have shown that
\begin{lemma}\label{lemma3.4}
\[
W_{n_0}(k,l)^\Gamma_0 /W_{n_0+1} (k,l)^\Gamma_0\subset  \bigoplus _{(\alpha,\beta,\gamma,\delta)\in I^{n_0}_{k,l}}(M_{2n_0}(\Gamma))_{(\alpha,\beta,\gamma,\delta)}  =M_{2n_0}(\Gamma)^{\oplus |I_{k,l}^{n_0}|},
\]
where $(M_{2n_0}(\Gamma))_{(\alpha,\beta,\gamma,\delta)}=M_{2n_0}(\Gamma)$ is the space of modular forms of weight $2n_0$ for $\Gamma$,  and $I^{n_0}_{k,l}$ is the collection of four-tuples of partitions $(\lambda,\mu,\nu,\chi)$ with the part $n_0$, conformal weight $k$ and fermionic charge $l$.
\end{lemma}

\noindent Let $\alpha_{n_0}:
  W_{n_0}(k,l)^\Gamma_0  \longrightarrow    M_{2n_0}(\Gamma)^{\oplus |I_{k,l}^{n_0}|},$ be the projection map
 \[ 
 \sum_{n\geq n_0} \sum_{(\lambda,\mu,\nu,\chi)\in I^{n}_{k,l} } a_{-\lambda}\phi_{-\mu}  \psi_{-
 \nu} b_{-\chi}f_{\lambda,\mu,\nu,\chi}  \longmapsto      \{f_{\lambda,\mu,\nu,\chi} \}_{(\lambda,\mu,\nu,\chi)\in I^{n_0}_{k,l}}.
  \]
 
 \noindent Our first main result is
 
  \begin{theorem} \label{theorem3.5}
 We have the short exact sequence:
\begin{equation}\label{3.20}
 0\longrightarrow W_{n_0+1}(k,l)^\Gamma_0 \longrightarrow W_{n_0}(k,l)^\Gamma_0 \overset{\alpha_{n_0}}{\longrightarrow} M_{2n_0}(\Gamma)^{\oplus |I_{k,l}^{n_0}|}\longrightarrow 0. 
 \end{equation}
 \end{theorem}
The proof of Theorem \ref{theorem3.5} will be given in Section \ref{Section 4}.
As a corollary of Lemma \ref{lemma3.4}, we have
   \begin{prop} 
 For any congruence subgroup $\Gamma$, we have
 \[
 \dim \Omega^{ch}(\mathbb H,\Gamma)_{k,l}<\infty, \;\;\text{ for any } k\geq 0, l\in \mathbb Z,
 \]
 where $\Omega^{ch}(\mathbb H,\Gamma)_{k,l}$ denotes the conformal weight $k$ and fermionic charge $l$ subspace of $\Omega^{ch}(\mathbb H, \Gamma)$.
 \end{prop}
 
\noindent {\it Proof: }
By Lemma \ref{lemma3.4}, the dimension of $W_{n_0} (k,l)^\Gamma_0$ is bounded by  $\dim W_{n_0+1}(k,l)^\Gamma_0 +|I_{k,l}^{n_0}|\dim M_{2n_0}(\Gamma)$.
From the argument below (\ref{3.18}), for any four-tuple $(\lambda,\mu,\nu,\chi) $ with fixed conformal weight $k$ and fermionic charge $l$, $ |\lambda|+|\mu|+|\nu|+|\chi| $ is bounded from above, say bounded by $N$. Then  
\[
p(\lambda)-p(\mu)+p(\nu)-p(\chi) = -p(\mu)-p(\chi)  \geq -N ,\] 
which implies $W_{n}(k,l)^\Gamma _0=0$ for $n \geq N+1$.
Hence we prove the result by induction.\qed

\section{Lifting of Modular Forms to $\Omega^{ch}(\mathbb H,\Gamma)$} \label{Section 4}

In this section, we will study the lifting of modular forms to $\Omega^{ch}(\mathbb H,\Gamma)$ under the map $\alpha_{n_0}$ in (\ref{3.20}). 
We will prove Theorem \ref{theorem3.5} by giving the concrete lifting formulas, whose holomorphic functions consist of derivatives of modular forms when $n_0 \geq 1$, and derivatives of Eisenstein series together with the constant modular form when $n_0=0$.

\subsection{Lifting of Non-constant Modular Forms}\label{Section 4.1}

Let $\mathfrak g$ be the direct sum of Heisenberg Lie algebra and Clifford Lie superalgebra with the centers identified and $N=1$ in Section \ref{Section 2}. So $\mathfrak g$ is spanned by the odd elements $\phi_{m},\psi_{m}$, the even elements $a_{m},b_{m}$ and the center $C$ as follows
\[
\mathfrak g=Span_\mathbb C\{ a_{m},b_{m},\phi_{m},\psi_{m}\;|  m\in \mathbb Z \} \oplus \mathbb C C,
\]with the Lie bracket
\begin{equation}\label{4.1}
[a_{m},b_{n}]=  \delta_{m+n,0} C,\;\;  [\phi_{m},\psi_{n}]=\delta_{m+n,0}C
\end{equation}
Let $\mathcal U$ be the quotient of the universal enveloping algebra of $\mathfrak g$ by the ideal generated by $C=1$. 
Then $\mathcal U =\oplus _{n\in\mathbb Z}\mathcal U_n$, where $\mathcal U_n$ denotes the conformal weight $n$ subspace of $\mathcal U$. 
Let $\mathcal U_{n}^k =\sum_{i\leq k} \mathcal U_{n-i}\mathcal U_{i}$, then 
\[
\mathcal U_{n}^k \subset \mathcal U_n^{k+1}.
\]
Then the fundamental system $\{ \mathcal U_n^k\}_{k\in\mathbb Z}$
 with the condition that 
 \[\cap_{n\in \mathbb Z} \mathcal U_n^k =\{0\}, \;\;\; \cup_{n\in \mathbb Z} \mathcal U_n^k=\mathcal U_n,
 \]
  gives a linear topology on ${\mathcal U}_n$.
Denote by $\bar{\mathcal U}_n$ the completion of $\mathcal U_n$ with respect to the topology. 
And $\bar{\mathcal U}_n$ has a fundamental system $\{ \bar{\mathcal U}_n^k\}_{k\in \mathbb Z}$ of neighborhoods of $0$. Then the direct sum $\bar {\mathcal U}:= \oplus _{n\in \mathbb Z }   \bar {\mathcal U}_n$ is a complete topological ring (see similar constructions in \cite{FZ}).
Note that $\bar{\mathcal U} $ acts on $\Omega^{ch}(\mathbb H)$.

Define 
\[
R:= \mathbb C[\tau][(m \tau+n)^{-1}|   (m,n)\in \mathbb R^2\backslash\{(0,0)\} ].
\]
Then we can easily show that $R$ is closed under the linear fractional transformation
\[
\tau \mapsto \dfrac{\alpha \tau+\beta}{\gamma \tau +\delta}.
\]
 Since the conformal weight of the operators $a_n,\phi_n,\psi_n,b_n$ equals $-n$, and the group action preserves the conformal weight, we can check that all but finitely many terms in (\ref{3.8})-(\ref{3.11}) are contained in $\bar{\mathcal U}_{-n}^N  \otimes _\mathbb C R$, for arbitrary $N\in \mathbb Z$. 
Thus the adjoint action of $g$ on $a_n,\phi_n,\psi_n$ and $b_n$ for $n\in\mathbb Z$ are contained in $ \bar{\mathcal U} \otimes _\mathbb C R$ and hence the Lie group $SL(2,\mathbb R)$ acts on $ \bar{\mathcal U}\otimes _\mathbb C R$.  
Let $K$ be the left ideal in $\bar{\mathcal U}$ generated by the elements $a_{n},b_{n},\phi_{n}, \psi_{m}$ for $n\geq 1$ and $m\geq 0$. 
Then $K\otimes _\mathbb C R$ is preserved under the $SL(2,\mathbb R)$-action by (\ref{3.8})-(\ref{3.11}). 
Therefore $(\bar{\mathcal U}/K) \otimes _\mathbb C R$ has an $SL(2,\mathbb R)$-module structure. For $f\in \mathcal O(\mathbb H) \subset \Omega^{ch}(\mathbb H)$, because $a_nf=\phi_n f=b_nf=0$, $\psi_m f=0$ for $n\geq 1,m\geq 0$ , so $Kf =0$. 
 Therefore we have a map 
 \begin{align*}
 ((\bar{\mathcal U}/K) \otimes _\mathbb C R)\times \mathcal O(\mathbb H) & \longrightarrow \Omega^{ch}(\mathbb H)\\
 ( (u+K)\otimes r) f&\longmapsto urf,
 \end{align*}
  which is $SL(2,\mathbb R)$-equivariant in the sense that 
\[
\pi(g)  Ar  f =(\pi(g)(  A \otimes r)\pi(g)^{-1}) \pi(g) f, \;\;\;\text{ for any } r\in R, A\in \bar{\mathcal U}/K,  g\in SL(2,\mathbb R). 
\]
According to PBW theorem, $\bar{\mathcal U}/K$ has the following basis
\begin{equation}\label{4.2}
a_{-\lambda}\phi_{-\mu} \psi_{-\nu} b_{-\chi}a_0^kb_0^l \;\; \; k,l\in \mathbb Z_{\geq 0}.
\end{equation}

  \
  
 For an arbitrary four-tuple $(\lambda_0,\mu_0,\nu_0,\chi_0)$ of conformal weight $k$ and fermionic charge $l$, such that $-p(\lambda_0)+p(\mu_0)-p(\nu_0)+p(\chi_0)=n_0 \geq 1$, 
 we will first study the lifting of a nonconstant modular form $f$ of weight $2n_0$ to $W_{n_0}(k,l)^\Gamma_0$.
The idea is to find an operator $A\in \bar{\mathcal U}/K$, such that $Af \in W_{n_0}(k,l)^\Gamma_0$ with the leading term $a_{-\lambda_0}\phi_{-\mu_0}\psi_{-\nu_0}b_{-\chi_0}f$.

 Then the $\Gamma$-invariance of $Af$ implies
 \[
 (\pi(g)A \pi(g)^{-1}) f(gb)=(\pi(g)A\pi(g)^{-1})(\gamma b+\delta)^{2n_0}f(b),\;\;\text{ for } g\in \Gamma.
 \] 
 So it suffices to find solutions of an operator $A$ such that 
 \begin{equation}\label{4.3}
 \pi(g)A\pi(g)^{-1}=A(\gamma b+\delta)^{-2n_0}, \;\;\; \text{ for any } g\in  SL(2,\mathbb R).
 \end{equation}
  Considering the infinitesimal action, (\ref{4.3}) is equivalent to the following system
 \begin{align} 
 \label{4.4} E_{(0)}. A&=0,\\
\label{4.5} H_{(0)}.A &=-2n_0 A,\\
\label{4.6} F_{(0)}.A &=2n_0A b_0,
 \end{align}
 where for $x\in \mathfrak{sl}_2$, we denote by $x_{(0)}.$ the infinitesimal adjoint action of $x$ on $\bar{\mathcal U}/K$, which is given by:  
 \begin{equation}\label{4.7}
x_{(0)}. B= x_{(0)}B-Bx_{(0)},\;\;\; \text{ for any } B\in \bar{\mathcal U}/K.
 \end{equation}
 
\noindent The adjoint action of $E_{(0)}$ and $H_{(0)}$ on the operator $a_{-\lambda} \phi_{-\mu} \psi_{-\nu}b_{-\chi} a_0^k b_0^l$ can be calculated easily, 
\begin{gather}\label{4.8}
E_{(0)}. a_{-\lambda}\phi_{-\mu}\psi_{-\nu} b_{-\chi} a_0^k b_0^l  =-l a_{-\lambda}\phi_{-\mu}\psi_{-\nu} b_{-\chi} a_0^k b_0^{l-1}, \\
\label{4.9}H_{(0)}. a_{-\lambda}\phi_{-\mu}\psi_{-\nu} b_{-\chi} a_0^k b_0^l  =2(p(\lambda)-p(\mu)+p(\nu)-p(\chi)+k-l) a_{-\lambda}\phi_{-\mu}\psi_{-\nu} b_{-\chi} a_0^k b_0^l .
\end{gather}
But the $F_{(0)}.$-action is quite complicated. 
Instead we will first consider the induced action on the graded algebra 
\[
Gr \,\Omega^{ch}(\mathbb H)=\oplus_{n\in \mathbb Z} W_n/W_{n+1} ,
\]
which is isomorphic to the tensor product $\Omega_1'\otimes _{\mathbb C} O(\mathbb H)$. And for $a_{-\lambda} \phi_{-\mu} \psi_{-\nu} b_{-\chi} f \in Gr \,\Omega^{ch}(\mathbb H)$, the induced $\mathfrak{sl}_2$-action is given as follows
\begin{subequations}
\begin{align}
\label{4.10a} E_{(0)} a_{-\lambda} \phi_{-\mu} \psi_{-\nu} b_{-\chi} f=& a_{-\lambda} \phi_{-\mu} \psi_{-\nu} b_{-\chi} (-f'),\\
\label{4.10b}H_{(0)} a_{-\lambda} \phi_{-\mu} \psi_{-\nu} b_{-\chi} f=& a_{-\lambda} \phi_{-\mu} \psi_{-\nu} b_{-\chi}(-2nf-2bf'),\\
\label{4.10c} F_{(0)} a_{-\lambda} \phi_{-\mu} \psi_{-\nu} b_{-\chi}f =& a_{-\lambda} \phi_{-\mu} \psi_{-\nu} b_{-\chi} (2nbf+b^2f'),
\end{align}
\end{subequations}
where $n=-p(\lambda)+p(\mu)-p(\nu)+p(\chi)$. Notice that from (\ref{4.8}) and (\ref{4.9}) we can see that the action formulas of $E_{(0)}$ and $H_{(0)}$ on $\Omega^{ch}(\mathbb H)$ coincide with (\ref{4.10a}) and (\ref{4.10b}), while the action formula of $F_{(0)}$ on $\Omega^{ch}(\mathbb H)$ involves more terms comparing to (\ref{4.10c}).

Let $C$ be the Casimir operator defined by
\begin{align}
\nonumber C=& E_{(0)} F_{(0)}+F_{(0)}E_{(0)} +\frac{1}{2} H_{(0)}^2\\
=&2F_{(0)}E_{(0)} +H_{(0)} +\frac{1}{2} H_{(0)}^2.\label{4.11}
\end{align}
From (\ref{4.10a})-(\ref{4.10c}) and (\ref{4.11}), we can compute the action formula of $C$ on $Gr \,\Omega^{ch}(\mathbb H)$ explicitly, 
\begin{equation}
C a_{-\lambda} \phi_{-\mu} \psi_{-\nu} b_{-\chi} f= a_{-\lambda} \phi_{-\mu} \psi_{-\nu} b_{-\chi} (2n(n-1)f).
\end{equation}
Hence $C$ acts semisimply on $Gr \,\Omega^{ch}(\mathbb H)$, and it acts as $2n(n-1)Id$ on $Gr \,\Omega^{ch}(\mathbb H)_n=W_n/W_{n-1}$. 

\begin{lemma}\label{lemma4.1}
The system (\ref{4.4})-(\ref{4.6}) is equivalent to the same system with (\ref{4.6}) replaced by the following equation
\begin{equation}\label{4.13}
C A f= 2n_0 (n_0-1) Af,\;\;\; \text{ for any } f\in \mathcal O(\mathbb H).
\end{equation}
\end{lemma}
\noindent {\it Proof:} 
By (\ref{4.5}), we have
\[
H_{(0)} Af=(H_{(0)}.A)f+AH_{(0)}.f= -2n_0Af+ 2Abf'.
\]
Similarly we may calculate
\begin{equation}\label{4.14}
(H_{(0)}+\frac{1}{2}H_{(0)}^2) Af =2n_0(n_0-1) Af  +4n_0Abf'+2Ab^2 f''.
\end{equation}
By (\ref{4.4}), 
\begin{equation}\label{4.15}
F_{(0)} E_{(0)} Af=-(F_{(0)}.A)f' -A F_{(0)}f'=  -(F_{(0)}.A)f' -A(b^2f'').
\end{equation}
Adding twice of equation (\ref{4.15}) to (\ref{4.14}) and using (\ref{4.11}), we have
\begin{equation} \label{4.16}
CAf=2n_0(n_0-1)Af +4n_0 Abf' -2 (F_{(0)}. A) f'.
\end{equation}
Therefore (\ref{4.13}) is equivalent to 
\[
(F_{(0)}.A) f'= 2n_0 Ab f',\;\;\;\text{ for all }  f\in \mathcal O(\mathbb H),
\]
which is equivalent to (\ref{4.6}).
\qed

Let 
\begin{equation} \label{4.17}
D: =F_{(0)}. + b H_{(0)}.
\end{equation}
 be an operator on $\bar{\mathcal U}/K$.
Observe that  $D$ is locally nilpotent. 
 For any $a_{-\lambda}\phi_{-\mu}\psi_{-\nu}b_{-\chi}f\in \Omega^{ch}(\mathbb H)$ with $n=-p(\lambda)+p(\mu)-p(\nu)+p(\chi)$,
\begin{align}\nonumber
F_{(0)}  a_{-\lambda}\phi_{-\mu}\psi_{-\nu}b_{-\chi}f=& (F_{(0)}.a_{-\lambda}\phi_{-\mu}\psi_{-\nu}b_{-\chi}) f+a_{-\lambda}\phi_{-\mu}\psi_{-\nu}b_{-\chi} F_{(0)} f\\
\label{4.18}=& D(a_{-\lambda}\phi_{-\mu}\psi_{-\nu}b_{-\chi}) f+ a_{-\lambda}\phi_{-\mu}\psi_{-\nu}b_{-\chi} (2n bf +b^2 f'),
\end{align}
which implies that the action formula of $F_{(0)}$ on $\Omega^{ch}(\mathbb H)$ differs from (\ref{4.10c}) by the action $D$ on the corresponding vertex operators.

Now since $A: =\sum_{n=0}^\infty c_n D^{n}(a_{-\lambda_0}\phi_{-\mu_0} \psi_{\nu_0} b_{-\chi_0} ) a_0^n$ automatically satisfies (\ref{4.4})-(\ref{4.5}), it suffices to find a family of constants $c_0=1,c_1,c_2,\cdots$, such that 
\begin{equation}\label{4.19}
C \sum_{n=0}^\infty c_n D^n(a_{-\lambda_0}\phi_{-\mu_0} \psi_{-\nu_0} b_{-\chi_0}) f^{(n)}=2n_0(n_0-1) \sum_{n=0}^\infty c_n D^n (a_{-\lambda_0}\phi_{-\mu_0} \psi_{-\nu_0} b_{-\chi_0}) f^{(n)}.
\end{equation}
Similar to the calculation of (\ref{4.16}), the left side of (\ref{4.19}) equals 
\[
\sum_{n=1}^\infty (2(n_0+n)(n_0+n-1)c_n-2c_{n-1}) D^{n}(a_{-\lambda_0}\phi_{-\mu_0} \psi_{-\nu_0} b_{-\chi_0}) f^{(n)} +2n_0(n_0-1) a_{-\lambda_0}\phi_{-\mu_0} \psi_{-\nu_0} b_{-\chi_0}f
\]
Hence we have recursive relations
$2(n_0+n)(n_0+n-1)c_n-2c_{n-1}=2n_0(n_0-1) c_{n}$ for $n\geq 1$, which may be rewritten as follows
\begin{equation}
c_n=\frac{1}{n(n+2n_0-1)}c_{n-1}, \text{ for } n\geq 1,
\end{equation}
and the first term $c_0=1$. Thus 
\[
c_{n}=\frac{(2n_0-1)!}{n!(n+2n_0-1)!}, \text{ for } n\geq 0.
\]
 So we have proved the lifting theorem for non-constant modular forms.

\subsection{Lifting of Constant Modular Form}\label{Section 4.2}

Recall that the Eisenstein series 
\[
E_2(\tau)=1-24\sum_{n=1}^\infty \sigma(n) q^n\;\;\; q=e^{2\pi i \tau}
\]
is a quasi-modular form of weight $2$, with the transformation property(cf.\cite{Z1}p.19)
\begin{equation} \label{4.21}
E_2(g\tau) =(\gamma \tau+\delta)^2 E_2(\tau)-\frac{6i}{\pi} \gamma(\gamma \tau +\delta),\;\;\; \text{ for any } g=\begin{pmatrix} \alpha & \beta \\ \gamma & \delta \end{pmatrix}\in SL(2,\mathbb Z).
\end{equation}

Define $E(b):=\frac{\pi i}{6} E_2(b) $, then (\ref{4.21}) is equivalent to
\begin{equation} \label{4.22}
(\gamma\tau +\delta)^{-2}E(g\tau)=E(b) +  \gamma (\gamma \tau +\delta)^{-1}.
\end{equation}

Since modular forms of weight $0$ are constant functions, it suffices to study the lifting of the constant function $ f(b)\equiv  1$. 
Now we will show that there exists an operator $A\in \bar{U}/K$, such that 
\[
\pi(g)  (a_{-\lambda_0} \phi_{-\mu_0}\psi_{-\nu_0} b_{-\chi_0} +AE(b)     )=a_{-\lambda_0}\phi_{-\mu_0}\psi_{-\nu_0}b_{-\chi_0} +AE(b),\;\;\;\;\text{ for any } g\in \Gamma.
\]
Notice that 
\[\pi(g) AE(b)=(\pi(g)A\pi(g)^{-1}) E(gb)=(\pi(g)A\pi(g)^{-1})((\gamma b+\delta)^2E(b)+\gamma(\gamma  b+\delta)).\]
It is natural to consider solutions of $A$, such that for $g\in SL(2,\mathbb R)$
\begin{align}
 \label{4.23}
\pi(g)A\pi(g)^{-1}&=A(\gamma b+\delta)^{-2},\\   
\label{4.24} -A \gamma (\gamma b+\delta)^{-1} &=(\pi(g)-I) a_{-\lambda_0 }\phi_{-\mu_0}\psi_{-\nu_0}b_{-\chi_0}.
\end{align}

Similar to (\ref{4.3}), the equation (\ref{4.23}) can be replaced by the version of infinitesimal action, namely 
\begin{align}
\label{4.25}E_{(0)}. A&=0,\\
\label{4.26}H_{(0)}. A&=-2A,\\
\label{4.27}F_{(0)}. A&=2Ab.
\end{align}

From Section \ref{Section 4.1}, solutions for operators $A\in \bar{U}/K$ satisfying (\ref{4.23}) exist.  For any such fixed operator $A$,
 we define a twisted $SL(2,\mathbb R)$-action on $\Omega^{ch}(\mathbb H)$ as follows,
 \begin{equation}\label{4.28}
\pi_A(g) v=\pi(g)v +A\gamma (\gamma b+\delta)^{-1},
\end{equation}
where $v\in \Omega^{ch}(\mathbb H)$ and $g\in SL(2,\mathbb R)$.
Obviously $\pi_A(I_2)$ acts as the identity operator on $\Omega^{ch}(\mathbb H)$, where $I_2$ is the identity matrix. Take arbitrary $g_i= \begin{pmatrix} \alpha_i & \beta_i \\ \gamma_i & \delta_i \end{pmatrix} \in SL(2,\mathbb R) $, for $i=1,2$, and denote by $g_1g_2=\begin{pmatrix} \alpha & \beta \\ \gamma & \delta \end{pmatrix}$. Then the compatibility condition can be derived as below
\begin{align*}
\pi_A(g_2)\pi_A(g_1) v=&\pi_A(g_2) (\pi(g_1)v+A\gamma_1(\gamma_1 b+\delta_1)^{-1}) \\
=&\pi(g_2)\pi(g_1) v+\pi(g_2)A \pi(g_2)^{-1} \pi(g_2)\gamma_1 (\gamma_1 b+\delta_1)^{-1}+ A\gamma_2 (\gamma_2 b+\delta_2)^{-1}\\
=&\pi(g_1g_2) v+A \gamma_1 (\gamma b+\delta)^{-1} (\gamma_2 b+\delta_2)^{-1} +A\gamma_2 (\gamma_2 b+\delta_2)^{-1}\\
=&\pi(g_1g_2)v+ A\gamma(\gamma b+\delta)^{-1}=\pi_A(g_1g_2) v,
\end{align*}
where we use (\ref{4.23}) in the third equation. Hence $\pi_A$ gives a well-defined right $SL(2,\mathbb R)$-action on $\Omega^{ch}(\mathbb H)$. Then (\ref{4.24}) can be rewritten as
\begin{equation} \label{4.29}
\pi_A(g) a_{-\lambda_0} \phi_{-\mu_0}\psi_{-\nu_0} b_{-\chi_0} = a_{-\lambda_0} \phi_{-\mu_0}\psi_{-\nu_0} b_{-\chi_0}.
\end{equation}
Now we will consider the twisted infinitesimal action of (\ref{4.28}) on $\Omega^{ch}(\mathbb H)$. Since the twisted action $\pi_A$ coincides with $\pi$ when restricting to the Borel subgroup consisting of the upper triangular matrices in $SL(2,\mathbb R)$, so does the twisted infinitesimal action with the original infinitesimal action of the corresponding Borel subalgebra. Therefore we have
\[ \pi_A  \begin{pmatrix} 0 & 1 \\ 0 & 0 \end{pmatrix} =E_{(0)} ,\;\;\; \pi_A \begin{pmatrix} 1 & 0 \\ 0 & -1 \end{pmatrix}  =H_{(0)}. \]
 And by a simple calculation, the twisted infinitesimal action of $\begin{pmatrix} 0 & 0\\ 1 & 0 \end{pmatrix}$ differs from $F_{(0)}$ by a translation of the vector $-A1$, namely for any $v\in \Omega^{ch}(\mathbb H)$, 
\[\pi_A \begin{pmatrix} 0 & 0\\ 1 & 0 \end{pmatrix} v=F_{(0)}v-A1.\]
Since $E_{(0)}$ and $H_{(0)}$ kill $a_{-\lambda_0 }\phi_{-\mu_0}\psi_{-\nu_0}b_{-\chi_0}$ automatically, the only nontrivial equation for the twisted infinitesimal version of (\ref{4.29}) is
\begin{equation}
\label{4.30}F_{(0)} a_{-\lambda_0 }\phi_{-\mu_0}\psi_{-\nu_0}b_{-\chi_0}=A1.
\end{equation}
 
 Hence we have proved the following lemma.

\begin{lemma}\label{lemma4.2}
The conditions (\ref{4.23}) and (\ref{4.30}) imply (\ref{4.24}).
\end{lemma}

By Lemma \ref{lemma4.2}, we only need to find $A\in \bar{U}/K$, satisfying (\ref{4.25})-(\ref{4.27}), and (\ref{4.30}). 
\begin{lemma}
The system (\ref{4.25})-(\ref{4.27}) is equivalent to the same system with (\ref{4.27}) replaced by the following equation
\begin{equation}\label{4.31}
C A f= 0,\;\;\; \text{ for any } f\in \mathcal O(\mathbb H).
\end{equation}
\end{lemma}
The proof is an analogy to the proof of Lemma \ref{lemma4.1}.

Now since the operator
\begin{equation}\label{4.32}
A: =\sum_{n=1}^\infty d_n D^{n}(a_{-\lambda_0}\phi_{-\mu_0} \psi_{\nu_0} b_{-\chi_0} ) a_0^{n-1}
\end{equation}
 satisfies (\ref{4.25})-(\ref{4.26}), we still need to find a family of constant $d_1,d_2,\cdots$, such that 
\begin{equation}\label{4.33}
C \sum_{n=1}^\infty d_n D^n(a_{-\lambda_0}\phi_{-\mu_0} \psi_{-\nu_0} b_{-\chi_0}) f^{(n-1)}=0.
\end{equation}
Using (\ref{4.14}) and (\ref{4.15}), the left side of (\ref{4.33}) equals 
\[
\sum_{n=1}^\infty (2n(n-1)d_n-2d_{n-1}) D^{n}(a_{-\lambda_0}\phi_{-\mu_0} \psi_{-\nu_0} b_{-\chi_0}) f^{(n-1)} 
\]
Hence we have the recursive relations
\begin{equation}
d_n=\frac{1}{n(n-1)}d_{n-1}, \text{ for } n\geq 1.
\end{equation}
Applying (\ref{4.32}) to (\ref{4.30}), we have
\begin{equation}\label{4.35}
F_{(0)}a_{-\lambda_0}\phi_{-\mu_0} \psi_{-\nu_0} b_{-\chi_0}= d_1 D(a_{-\lambda_0}\phi_{-\mu_0} \psi_{-\nu_0} b_{-\chi_0}) 1.
\end{equation}
Thanks to (\ref{4.18}), the left side of (\ref{4.35}) equals
\[
D(a_{-\lambda_0}\phi_{-\mu_0} \psi_{-\nu_0} b_{-\chi_0}) 1.
\]
If $D(a_{-\lambda_0}\phi_{-\mu_0} \psi_{-\nu_0} b_{-\chi_0}) 1=0$, then (\ref{4.35}) trivially holds for arbitrary $d_1$. Otherwise, $d_1=1$. In either case, we take $d_1=1$,
and hence $d_{n}=\frac{1}{n!(n-1)!}$, for $n\geq 1$. So we have proved the lifting theorem for the non-constant modular forms.

\section{Structures of the Vertex Algebra $\Omega^{ch}(\mathbb H,\Gamma)$} \label{Section 5}
\subsection{Topological Vertex Algebra Structure}
For any nonconstant modular form $f\in M_{2k}(\Gamma)$,
and arbitrary operator $w=a_{-\lambda}\phi_{-\mu}\psi_{-\nu} b_{-\chi}$ with the part $k$, namely $-p(\lambda)+p(\mu)-p(\nu)+p(\chi)=k$, 
we define $L(w,f)$ to be the unique lifting of $f$ in $\Omega^{ch}(\mathbb H,\Gamma)$ with the leading term $wf$ as in Theorem \ref{theorem1.1}, namely
\begin{equation}
L(w, f):=\sum_{n=0}^\infty \frac{(2k-1)!}{n!(n+2k-1)!} D^{n}(w) f^{(n)},
\end{equation}
where the operator $D $ is defined as in (\ref{4.17}). When the part of $w=a_{-\lambda}\phi_{-\mu}\psi_{-\nu} b_{-\chi}$ is not equal to $k$, we define $L(w,f)$ to be $0$.
And for any constant modular form $c\neq 0$, we define
\begin{equation}
L(w,c):= cw+\sum_{n=1}^\infty \frac{c}{n!(n-1)!} D^{n}(w) E^{(n-1)}\in\Omega^{ch}(\mathbb H,\Gamma(1)),
\end{equation}
to be the lifting of $c$ in $\Omega^{ch}(\mathbb H,\Gamma)$ with the leading term $cw$ when the operator $w$ has the part $0$. Similarly we define $L(w, c)=0$, if the part of $w$ is nonzero. 
And then we extend the definition of the $L$ operator to the whole space $\Omega'_1\otimes M(\Gamma)$ by linearity, where $M(\Gamma)=\oplus_{n\geq 0} M_{2n} (\Gamma)$.

\begin{lemma} \label{lemma5.1}
Let $B_{2k}(\Gamma)\subset M_{2k}(\Gamma)$ be a linear basis of modular forms of weight $2k$, for $k\geq 0$, and $I_k$ be the collection of all four-tuples of partitions $(\lambda,\mu,\nu,\chi)$ with $-p(\lambda)+p(\mu)-p(\nu)+p(\chi)=k$. Then the following liftings
\begin{equation}\label{5.3}
\{ L(a_{-\lambda}\phi_{-\mu} \psi_{-\nu} b_{-\chi}, f) \;|\; (\lambda,\mu,\nu,\chi)\in I_{k}, f\in \mathcal B_{2k}(\Gamma) , \text{ for all } k\geq0 \}.
\end{equation}
 form a linear basis in $\Omega^{ch}(\mathbb H,\Gamma)$.
\end{lemma}

We can easily verify that the Virasoro element $\omega= b_{-1}a_{-1}+\phi_{-1} \psi_{-1}$ and $G=\psi_{-1} b_{-1}$ are still fixed by $SL(2,\mathbb R)$, which means $\omega$ and $G$ are contained in $\Omega^{ch}(\mathbb H,\Gamma)$ for any congruence subgroup $\Gamma$. However the even element $J=\phi \psi_{-1}$ and the odd element $Q=a_{-1} \phi$ are not contained in $\Omega^{ch}(\mathbb H,\Gamma)$ in general. As a first application of Theorem \ref{theorem1.1}, we consider the liftings of the constant modular form $1$ with the leading terms $J$ and $Q$.
Denote by $\tilde J=L(J,1)$, and $\tilde Q=L(Q,1)$.
Applying (\ref{1.3}),
we obtain (\ref{1.5}) and (\ref{1.6}).
We can easily check that the corresponding fields $\tilde J(z)$, $\tilde Q(z)$, together with the fields $L(z)$ and $G(z)$ satisfy the relations (\ref{2.14})-(\ref{2.17}), so they make $\Omega^{ch}(\mathbb H,\Gamma)$ a topological vertex algebra. 
And notice that as a derivation, $(b_{-1}E)_{(0)}$ acts trivially on the generators $a,b_{-1},\phi,\psi$ and $f\in \mathcal O(\mathbb H)$, hence it acts as a zero operator on $\Omega^{ch}(\mathbb H,\Gamma)$. So $\tilde J_{(0)}=J_{(0)}$.
Therefore $\tilde J_{(0)}$ acts semisimply on $\Omega^{ch}(\mathbb H,\Gamma)$, and the eigenvalue is exactly the fermionic charge, namely
\[
\tilde J_{(0)} L(a_{-\lambda}\phi_{-\mu} \psi_{-\nu} b_{-\chi},f)=( p(\mu)-p(\nu))L(a_{-\lambda}\phi_{-\mu} \psi_{-\nu} b_{-\chi},f).
\]
 And since $\phi_{-1}E+\phi b_{-1} E' =T(\phi E )$, we also have $\tilde Q_{(0)}=Q_{(0)}$.  Let $d^{ch}=-\tilde Q_{(0)}$. Then the operator $d^{ch}$  increases the fermionic charge by one and
 $\Omega^{ch}(\mathbb H,\Gamma)$ equipped with the chiral de Rham differential $d^{ch}$, forms a complex. 
  Let $\Omega^{ch}(\mathbb H,\Gamma)^0 \subset \Omega^{ch}(\mathbb H, \Gamma)$ spanned by conformal weight zero part, namely
   \[
 \Omega^{ch}(\mathbb H,\Gamma)^0 =\mathbb C1\oplus Span_\mathbb C \{ L(\phi, f) \; | \; f\in M_2(\Gamma)\},
 \]
  and we let $d'$ be the restriction of $d^{ch}$ to $\Omega^{ch}(\mathbb H,\Gamma)^0$. Then $ \Omega^{ch}(\mathbb H,\Gamma)^0$ 
  equipped with $d'$ forms a subcomplex of $ \Omega^{ch}(\mathbb H,\Gamma)$.
 According to the relation
  \[
 [ \tilde Q_{(0)}, G_{(1)} ]=L_{0},
 \]
 we have the following lemma.
 \begin{lemma}
The following embedding 
\[
(\Omega^{ch}(\mathbb H,\Gamma)^0,d') \longrightarrow (\Omega^{ch}(\mathbb H,\Gamma), d^{ch})
\]
is a quasi-isomorphism.
\end{lemma}

From the above lemma, we can compute the cohomology group of $\Omega^{ch}(\mathbb H,\Gamma)$ as follows
\begin{gather*}
H^0(\Omega^{ch}(\mathbb H,\Gamma))=\mathbb C1, \;\; H^{1}(\Omega^{ch}(\mathbb H,\Gamma))=M_2(\Gamma), \\
H^n(\Omega^{ch}(\mathbb H,\Gamma))=0, \text{ for } n\geq 2.
\end{gather*}

\
\

\subsection{Character Formula}

Now we will derive the character formula of $\Omega^{ch}(\mathbb H,\Gamma)$, which is the formal power series of variable $q$ defined by $\sum_{n=0}^\infty \dim \Omega^{ch}(\mathbb H,\Gamma)_n q^n$, namely
$tr\, q^{L_0}$.

We first consider the trace $tr\, t^{H_{(0)}}q^{L_0}$ of the vertex subalgebra $\Omega_1'$. According to (\ref{4.9}), $tr\, t^{H_{(0)}}q^{L_0}$ equals
\begin{equation}\label{5.4}
\prod _{n=1}^{\infty} \dfrac{1}{1-t^2q^n} \prod_{n=1}^\infty \dfrac{1}{1-t^{-2}q^n} \prod_{n=1}^\infty(1+t^2q^n) \prod_{n=0}^\infty (1+t^{-2} q^n)=\sum_{n=0}^\infty \sum_{m=-\infty}^{\infty} c(m,n) t^{m}q^n
\end{equation} 
It is clear that $c(m,n)$ is the number of four-tuples of partitions $(\lambda,\mu,\nu,\chi)$, satisfying the relations that $|\lambda|+|\mu|+|\nu|+|\chi|-p(\mu)=n$ and $2(p(\lambda)-p(\mu)+p(\nu)-p(\chi))=m$. 
Notice that $\dim M_{m}(\Gamma)$ equals zero when $m< 0$. By Lemma \ref{lemma5.1}, the character of $\Omega^{ch}(\mathbb H,\Gamma)$ equals
\begin{equation}\label{5.5}
\sum_{n=0}^\infty \sum_{m=-\infty}^\infty c(2m,n) \dim M_{2m}(\Gamma) q^n.
\end{equation}

\noindent  {\it Proof of Theorem \ref{theorem5.3.3}:}
The left side of (\ref{5.4}) equals
\begin{gather*}
(1+t^{-2})\sum_{s_1,s_2,s_3,s_4\geq 0} \sum_{n_1,n_2,n_3,n_4\geq 0} p_{s_1}(n_1)p_{s_2}(n_2) p'_{s_3}(n_3) p'_{s_4}(n_4) t^{2s_1-2s_2+2s_3-2s_4} q^{n_1+n_2+n_3+n_4}\\
=(1+t^{-2}) \sum_{s_1,s_2,s_3,s_4\geq 0} \left(q^{s_1} \prod_{i=1}^{s_1} \dfrac{1}{1-q^i} \right) \left(q^{s_2} \prod_{j=1}^{s_2} \dfrac{1}{1-q^j}\right) \left(q^{\frac{1}{2}s_3(s_3+1)}  \prod_{k=1}^{s_3} \dfrac{1}{1-q^k}    \right)  \left( q^{\frac{1}{2}s_4(s_4+1)} \prod_{l=1}^{s_4} \dfrac{1}{1-q^l}  \right) \\
\cdot t^{2s_1-2s_2+2s_3-2s_4}\\
=\sum_{s_1,s_2,s_3,s_4\geq 0} q^{s_1+s_2+\frac{1}{2}s_3(s_3+1)+\frac{1}{2}s_4(s_4-1)} \prod_{i=1}^{s_1} \dfrac{1}{1-q^i}  \prod_{j=1}^{s_2} \dfrac{1}{1-q^j} \prod_{k=1}^{s_3} \dfrac{1}{1-q^k}  \prod_{l=1}^{s_4} \dfrac{1}{1-q^l} \cdot t^{2s_1-2s_2+2s_3-2s_4},
 \end{gather*}
where $p_k(n)$ is the number of partitions of $n$ into exactly $k$ parts, $p'_k(n)$ is the number of partitions of $n$ into distinct $k$ parts, and in the second equality, we use the generating function of $p_k(n)$ and $p'_k(n)$.
Recall that a partition $\lambda$ has $k$ parts if and only if its conjugate partition $\lambda'$ has the largest part $k$, 
where the conjugate partition $\lambda'$ is the partition whose Young diagram is obtained by interchanging rows and columns of $\lambda$.
So the generating function for partition with part $k$, is
\[
\sum_{n\geq 0}p_k(n)x^n=x^k\prod_{i=1}^k  \frac{1}{1-x^i}.
\]
And since $p'_k(n)=\sum_{i=0}^k p_i \left(n-\frac{k(k+1)}{2}\right)$, the generating function for partition with distinct $k$ parts is
\[
\sum_{n\geq 0} p'_k(n)x^n = \sum_{n\geq 0} \sum_{i=0}^k p_i\left(n-\frac{k(k+1)}{2} \right) x^{n-\frac{k(k+1)}{2}}x^{\frac{k(k+1)}{2}}=x^{\frac{k(k+1)}{2}} \prod_{i=1}^k \dfrac{1}{1-x^i}.
\]

Hence (\ref{5.5}) equals
\begin{equation}\label{5.6}
\sum_{s_1,s_2,s_3,s_4\geq 0}^\infty  \dim M_{-2s_1+2s_2-2s_3+2s_4} (\Gamma) q^{s_1+s_2+\frac{1}{2}s_3(s_3+1)  +\frac{1}{2}s_4(s_4-1) } \prod_{i=1}^{s_1} \dfrac{1}{1-q^i}\prod_{j=1}^{s_2} \dfrac{1}{1-q^j} \prod_{k=1}^{s_3}\dfrac{1}{1-q^k} \prod_{l=1}^{s_4} \dfrac{1}{1-q^l}.
\end{equation}
(\ref{5.6}) is equivalent to (\ref{1.4}) by a change of variables $-s_1+s_2-s_3+s_4 \mapsto m$, $s_1+s_3\mapsto n$, $s_3 \mapsto u$ and $s_4 \mapsto v$.

\section{Rankin-Cohen Operators and Invariant Global Sections}  \label{Section 6}
In this section, we will briefly recall the definition of the Rankin-Cohen bracket of modular forms, which is a family of universal bilinear operations sending two modular forms to a modular form.  
As explained in  \cite{Z2} and \cite{CMZ},  the product of two $\Gamma$-invariant pseudodifferential operators is again a pseudodifferential operator whose components are scalar multiples of the Rankin-Cohen bracket, 
which implies a noncommutative multiplicative structure on the space of modular forms.
A similar idea can be applied to explore the connection between the invariant global sections and the Rankin-Cohen bracket. And we will see that the vertex operator structure will imply a modified Rankin-Cohen bracket on modular forms, 
which is nontrivial whenever constant modular forms involved.

\subsection{The Modified Rankin-Cohen bracket} \label{Section 6.1}
Let $\Gamma\subset SL(2,\mathbb Z)$ be a congruence subgroup, and $f\in M_{k}(\Gamma)$, $h\in M_{l}(\Gamma)$, then the $n$th Rankin-Cohen bracket is given by
\begin{equation}\label{6.1}
[f,h]_n=\frac{1}{(2\pi i)^n }\sum_{r+s=n} (-1)^r {n+k-1\choose s} {n+ l-1\choose r} f^{(r)}(\tau) h^{(s)}(\tau).
\end{equation}
 Therefore the graded vector space $M (\Gamma) =\oplus_{i\geq 0}M_{2i}(\Gamma)$ possesses an infinite family of bilinear operations $[\; ,\; ]_n: M_\ast \otimes M_{\ast}\mapsto M_{\ast+\ast+2n}$, with the $0$th bracket the usual multiplication. 
 The modularity of the Rankin-Cohen bracket was first proved by Cohen \cite{C}, involving the Cohen-Kuznetsov lifting from modular forms to the Jacobi-like forms (see in \cite{Z2} and \cite{CMZ}).
 Let $R$ be the ring of all holomorphic functions on $\mathbb H$ bounded by a power of $(|\tau|^2+1)/y$,
 \[
 R:=\{ f\in \mathcal O(\mathbb H) \;| \;  |  f(\tau)| \leq C (|\tau|^2+1)^l/y^l \;\;\text{ for some } l ,C>0\}.
 \]
 Then $\Gamma $ acts on $R$ via 
 \begin{equation}
f(\tau) \mapsto (\gamma b+\delta)^{-k}f(g\tau),
\end{equation}
 and
  the $\Gamma$-invariants $R^\Gamma $ is the space of modular forms of weight $k$ \cite{CMZ}. 
 A Jacobi-like form of weight $k$ for $\Gamma$ is defined to be a power series $\Phi(\tau,X)\in R[[X]]$, such that
  for any $g= \begin{pmatrix} \alpha & \beta \\ \gamma & \delta \end{pmatrix} \in \Gamma$, 
\begin{equation}
\Phi\left(\dfrac{\alpha \tau +\beta}{\gamma \tau +\delta}, \dfrac{X}{ (\gamma \tau +\delta)^2}\right) =(\gamma \tau +\delta)^k e^{\gamma  X/(\gamma \tau +\delta)} \Phi(\tau,X).
\end{equation}

We will denote by $J_k(\Gamma)$ the space of all Jacobi-like form of weight $k$ for the congruence subgroup $\Gamma$, and notice that the restriction of $\Phi$ to $X=0$, gives a modular form of weight $k$.
Now the Cohen-Kuznetsov lifting of a modular form $f$ is the formal power series
\begin{equation}\label{6.4}
\tilde f(\tau, X):=\sum_{n=0}^\infty \dfrac{f^{(n)}(\tau)} {{n! (n+k-1)!}}  X^n,
\end{equation}
whose Jacobi-like property can be found in \cite{Ku} and \cite{C}. 
Notice that when $k=2n_0$, the coefficient of $X^n$ in (\ref{6.4}) and the coefficients of $D^n(a_{-\lambda_0}\phi_{-\mu_0}\psi_{-\nu_0} b_{-\chi_0} )$ in  (\ref{1.2}) are exactly the same up to a fixed constant for all $n\geq 0$.

\

We assume that
 $w=a_{-\lambda_1}\phi_{-\mu_1} \psi_{-\nu_1}b_{-\chi_1},v=a_{-\lambda_2} \phi_{-\mu_2}\psi_{-\nu_2}b_{-\chi_2}\in \Omega'_1$ have the part $k$ and $l$ respectively, and $f_1,f_2$ are two modular forms for $\Gamma$ of weight $2k$ and $2l$ respectively. Since the $n$-th normal product of two liftings is still $\Gamma$-invariant, 
 there exists a sequence of modular forms associated to four-tuples of partitions, such that
\begin{equation}\label{6.5}
L(w, f_1)_{(n)} L(v,f_2)=\sum_{(\lambda,\mu,\nu,\chi)} L(a_{-\lambda} \phi_{-\mu} \psi_{-\nu} b_{-\chi},h_{\lambda,\mu,\nu,\chi}),
\end{equation}
where $h_{\lambda,\mu,\nu,\chi}$ is a modular form of weight $2s=-2p(\lambda)+2p(\mu)-2p(\nu)+2p(\chi)$. When $k,l\geq 1$, $(f_1,f_2)\mapsto h_{\lambda,\mu,\nu,\chi}$ is a map from $M_{2k}(\Gamma)\otimes M_{2l}(\Gamma)$ to $ M_{2s}(\Gamma)$, 
which can be written as a universal bilinear combination of products of derivatives of $f_1$ and $f_2$. Hence $h_{\lambda,\mu,\nu,\chi}$ is a multiple of the Rankin-Cohen bracket $[f_1,f_2]_{s-k-l}$. 
However when $k$ or $l$ equals $0$, the situation becomes more complicated. Actually if $f_1$ or $f_2$ is a constant modular form, the map $(f_1,f_2)\mapsto h_{\lambda,\mu,\nu,\chi}$ involves the derivatives of Eisenstein series $E_2$.
 This suggests to modify the Rankin-Cohen bracket when constant modular forms appear.

We first define the generalized Cohen-Kuznetsov lifting of the constant modular form $1$ to be
\begin{equation}\label{6.6}
\tilde{1}(\tau, X)=1+ \sum_{n=1}^\infty   \dfrac{E^{(n-1)}(\tau)}{n! (n-1)!} X^{n}.
\end{equation}
Formally the above formula can be viewed as (\ref{1.3}) with the operator $D^n(a_{-\lambda_0}\phi_{-\mu_0}\psi_{-\nu_0} b_{-\chi_0})$ replaced by $X^n$.

\begin{lemma}
The formula (\ref{6.6}) gives a Jacobi-like form of weight $0$ for the modular group $SL(2,\mathbb Z)$, i.e. $\tilde 1$ satisfies the transformation law
\begin{equation}\label{6.7}
\tilde 1\left(\dfrac{\alpha \tau +\beta}{\gamma \tau +\delta},  \dfrac{X}{(\gamma \tau +\delta)^2} \right)= e^{\gamma X/ (\gamma \tau +\delta)} \tilde{1}(\tau,X),\;\;\;\text{ for } g=\begin{pmatrix} \alpha & \beta \\ \gamma & \delta \end{pmatrix}\in SL(2,\mathbb Z).
\end{equation}
\end{lemma}
\noindent {\it Proof:} 
(\ref{6.7}) is equivalent to the following formula
\begin{equation} \label{6.8}
\dfrac{E^{(n)}(g \tau)}{n! (n+1)!} =\sum_{m=0}^n  \dfrac{ \gamma^{n-m} (\gamma \tau +\delta)^{m+n+2}}{(n-m)!}\dfrac{E^{(m)}(\tau)}{m!(m+1)!}+ \dfrac{\gamma^{n+1} }{(n+1)!} (\gamma \tau +\delta)^{n+1}, \;\;\;\;(n\geq 0).
\end{equation}
When $n=0$, (\ref{6.8}) is equivalent to (\ref{4.22}), and for general $n$, (\ref{6.8}) can be proved by induction.\qed

For any modular form $f\in M_{k}(\Gamma)$, with $k\geq 1$. The product of two Jacobi-like forms $\tilde 1$ and $\tilde f$ equals
\begin{equation}\label{6.9}
\tilde 1(\tau ,-X) \tilde f(\tau,X)=  \dfrac{f(\tau)}{(k-1)!}
+  \sum_{n\geq 0}  \left( \sum_{r+s=n}
 \dfrac{ (-1)^{r+1} E^{(r)}(\tau) f^{(s)}(\tau)}{r!s!(r+1)!(s+k-1)!}+ \dfrac{f^{(n+1)}}{(n+1)!(n+k)!}  \right)X^{n+1}, 
\end{equation}
whose coefficient of $X^n$ is a modular form of weight $k+2n$. 
If we make the substitution $r\to r-1$ to the coefficient of $X^n$, we have
\[
\dfrac{1}{(n+k-1)!(n-1)!}\left( \sum_{\substack{r+s=n\\ r\geq 1}} (-1)^r {n-1\choose s}{n+k-1\choose r} E^{(r-1)}(\tau) f^{(s)}(\tau) + \frac{f^{(n)}(\tau)}{n}\right) \in M_{k+2n}(\Gamma)
\]
Then we construct a family of linear maps for $n\geq 0,k>0$, 
\begin{align*}
M_k(\Gamma)& \longrightarrow M_{k+2n}(\Gamma),\\
f  &\longmapsto [1,f ]^\thicksim_n,
\end{align*}
where the modified bracket $[1,f]^\thicksim_n$ is defined by
\begin{equation}
[1, f ]^\thicksim _{n} :=\frac{1}{(2\pi i)^n}\sum_{\substack{r+s=n\\  1\leq r\leq n}  } (-1)^{r} {n-1 \choose s} {n+k-1 \choose r} E^{(r-1)}(\tau) f^{(s)}(\tau)+ \frac{1}{(2\pi i)^n}\dfrac{f^{(n)} (\tau)}{n},\;\;\;\text{ for } n> 0.
\end{equation}
We define $[1,f]^\thicksim _0:=f$, and require $[f, 1]^\thicksim_n= (-1)^{n}[1,f]^\thicksim_{n}$. As an example, we can calculate that
\[
[1,f]^\thicksim_1=\frac{1}{2\pi i} (-k E(\tau) f(\tau)+f'(\tau)) \in M_{k+2}(\Gamma).
\]

Now we will modify the $n$-th Rankin-Cohen bracket on two constant modular forms by multiplying the following two Jacobi-like forms, i.e.
\begin{equation}
\tilde{1}(\tau ,-X) \tilde 1(\tau,X) =1+\sum_{n\geq 0} \left( \sum_{r+s=n }\dfrac{ (-1)^{r+1} E^{(r)}(\tau) E^{(s)}(\tau) }{r! s!(r+1)!(s+1)!} + \dfrac{((-1)^{n} +1)E^{(n+1)}(\tau)   }{(n+1)!(n+2)!} \right)X^{n+2}.
\end{equation}
Substituting $r\to r-1$ and $s\to s-1$ to the coefficients of $X^n$, we have
\[
\dfrac{1}{((n-1)!)^2} \left( \sum_{\substack{r+s=n-2\\ 1\leq r \leq n-1}} (-1)^{r} {n-1\choose s+1 }{n-1\choose r+1}  E^{(r)}(\tau) E^{(s)}(\tau)+ \dfrac{(-1)^n+1}{n} E^{(n-1)}(\tau) \right) \in M_{2n}(\Gamma).
\]
This suggests to define the modified bracket $[1,1]^\thicksim_n$ as follows, for $n>0$
\begin{equation} \label{6.12}
[1,1]^\thicksim_n:= \frac{1}{(2\pi i)^n}\sum_{\substack{r+s=n\\ 1\leq r \leq n-1  }} (-1)^r {n-1\choose s}{n-1 \choose r} E^{(r-1)}(\tau) E^{(s-1)}(\tau)+\dfrac{(-1)^n+1}{(2\pi i)^n n} E^{(n-1)} (\tau),
\end{equation}
and $[1,1]^\thicksim_0:=1$. From (\ref{6.12}), it is easy to derive that
\[
[1,1]^\thicksim_{2k+1}=0,\;\;\; [1,1]^\thicksim_2=\frac{1}{(2\pi i)^2} (-E(\tau)^2 +E'(\tau)) \in M_{4}(\Gamma(1)).
\]
 
Now the $n$-th modified Rankin-Cohen bracket of modular forms are the bilinear operations: $M_\ast(\Gamma) \otimes M_\ast(\Gamma) \mapsto M_{\ast+\ast+2n}(\Gamma)$, defined as follows
\begin{equation}
[f,h]^\thicksim_n:=
\begin{cases}
[f,h]_n \;\;\;&\text{ if both}  f,h \text{ are nonconstants},\\
c[f,1]^\thicksim_n &\text{ if } f\in M_{k}(\Gamma), h=c , \text { for a constant $c$ and $k\geq 1$},  \\
d[1,h]^\thicksim_n & \text{ if } f=d, h\in M_{l}(\Gamma), \text{ for a constant $d$ and $l\geq 1$},\\
cd[1,1]^\thicksim_n &\text{ if $f=c,h=d$, for constants $c,d$.    }
\end{cases}
\end{equation}

Notice that if we formally view $E$ as the derivative of $1$, namely $E:=1'$ and hence $E^{(r)}=1^{(r+1)}$.  And we make the convention that $(-1)!:=1$, and redefine the combinatorial numbers  ${n-1\choose n}:=\frac{(n-1)!}{n! (-1)!}=\frac{1}{n}$ for $n\geq 1$. With this convention, ${-1 \choose 0}=\frac{(-1)!}{0! (-1)!}=1$ as usual. Then for $f\in M_{k}(\Gamma), h\in M_{l}(\Gamma)$, the modified Rankin-Cohen bracket can be rewritten as
\begin{equation}
[f,h]^\thicksim_n=\frac{1}{(2\pi i)^n}\sum_{r+s=n} (-1)^r { n+k-1\choose s}{n+l-1\choose r} f^{(r)}(\tau) h^{(s)}(\tau),
\end{equation}
which is exactly the same as (\ref{6.1}).

\begin{lemma}\label{lemma6.2}
$[\;,\;]^\thicksim_n$ is the unique universal combination of product of first $n$ derivatives of $1$(viewed as ``$-1$th derivative" of $E$) and a modular form of weight $k$ to a modular form of weight $k+2n$, namely, for any modular form $f\in M_{k}(\Gamma)$, if
\begin{equation} \label{6.15}
\sum_{r+s=n} c_{r,s} 1^{(r)} (\tau) f^{(s)}(\tau)\in M_{k+2n}(\Gamma),
\end{equation}
where $1^{(r)}:=E^{(r-1)}$ for $r\geq 0$, then (\ref{6.15}) equals a scalar multiple of $[1,f]^\thicksim_n$. 
\end{lemma}
\noindent {\it Proof:} 
Denote by $F$ the function in (\ref{6.15}). Then 
\begin{equation}\label{6.16}
F(g\tau)=(\gamma \tau+\delta)^{k+2n} F(\tau).
\end{equation}
The left hand side of (\ref{6.16}) equals
\begin{align*}
&\sum_{r+s=n} c_{r,s} E^{(r-1)}(g\tau) f^{(s)} (g\tau)\\
=& \sum_{r+s=n} \sum_{m=0}^{r} \sum_{t=0}^s c_{r,s} \dfrac{ (r-1)! r! s!(s+k-1)! \gamma ^{n-m-t}(\gamma \tau +\delta)^{k+n+m+t}} {(r-m)!(s-t)!(m-1)!m!   t!(t+k-1)! } E^{(m-1)}(\tau)f^{(t)}(\tau),
\end{align*}
where we use (\ref{6.8}), and transformation formula for modular form $f$ as follows
\[
\dfrac{f^{(n)}(g\tau)}{n! (n+k-1)!} =\sum_{m=0}^n \dfrac{\gamma^{n-m} (\gamma \tau +\delta)^{k +n+m} }{(n-m)!} \dfrac{f^{(m)}(\tau)}{m! (m+k-1)!},\;\;\; (n\geq 0),
\]
which can be proved by induction similar to (\ref{6.8}).
Now we compare the coefficient of the term $E^{(m-1)}(\tau)f^{(t)}(\tau)$ of both sides of (\ref{6.16}), we have
\[
c_{m,t} = \sum_{r+s=n}\sum_{ \substack{r\geq m\\s\geq t}}        c_{r,s} \dfrac{ (r-1)! r! s!(s+k-1)! \gamma ^{n-m-t}(\gamma \tau +\delta)^{m-n+t}} {(r-m)!(s-t)!(m-1)!m!   t!(t+k-1)! } .
\]
For any $r,s\geq 0$ such that $r+s=n$, take $m=r, t=s-1$. Then the above equation becomes
\[
c_{r,s} s(s+k-1) +c_{r+1,s-1} (r+1)r=0,
\]
which implies the coefficients $c_{r,s}$ are proportional to $(-1)^r {n-1\choose s}{n+k -1\choose r}$.\qed

From the discussion below (\ref{6.5}), the $n$-th normal product of two liftings $L(w,f)$ and $ L(v,h)$ is a linear combination of the liftings of the Rankin-Cohen bracket of $f$ and $h$, when $f$ and $h$ are nonconstant modular forms. 
When $f$ or $h$ is constant,  the modular form $h_{\lambda,\mu,\nu,\chi}$ in (\ref{6.5}) can be written as a combination of product of derivatives of $f$ and $h$ (here we note that $1^{(r)}=E^{(r-1)}$ as before), which 
by Lemma \ref{lemma6.2}, is a multiple of the modified Rankin-Cohen bracket of $f$ and $h$. And we have shown that
\begin{prop} \label{prop6.3}
For any modular forms $f_1\in M_{2k}(\Gamma),f_2\in M_{2l}(\Gamma)$, and any elements $w=a_{-\lambda_1}\phi_{-\mu_1} \psi_{-\nu_1} b_{-\chi_1}$, $v=a_{-\lambda_2} \phi_{-\mu_2} \psi_{-\nu_2} b_{-\chi_2}$ in $\Omega'_1$ of the part $k$ and $l$ respectively, we have
\begin{equation} 
L(w,f_1) _{(n)} L(v,f_2)=\sum_{(\lambda,\mu,\nu,\chi)} c_{\lambda,\mu,\nu,\chi}^nL(a_{-\lambda}\phi_{-\mu}\psi_{-\psi} b_{-\chi},  [f_1,f_2]^\thicksim_{-p(\lambda)+p(\mu)-p(\nu)+p(\chi)-k-l}),
\end{equation}
where $c_{\lambda,\mu,\nu,\chi}^n$ is a constant irrelevant to the choice of modular forms.
\end{prop}

Since $n$ is nonnegative, the weight of the modular form $[f_1,f_2]^\thicksim_n$ is always greater than or equal to $2k+2l$, namely the summation of the weights of $f_1$ and $f_2$.
Define
 \begin{equation}
L_n:=\text{Span}_\mathbb C \{ L(w,f) \in \Omega^{ch}(\mathbb H,\Gamma)\; |\; w\in I_k, f\in \mathcal B_{2k}(\Gamma) \text{ for } k\geq n\}.
\end{equation}
As a direct corollary of Proposition \ref{prop6.3}, $L_n$ is closed under the normal order product of elements in $\Omega^{ch}(\mathbb H,\Gamma)$. Namely we have 
\begin{cor}
$L_n$ is an ideal of vertex algebra $\Omega^{ch}(\mathbb H,\Gamma)$, for arbitrary $n\geq 0$
\end{cor}

 Let $\mathcal W\subset \Omega'_1$ be a vertex subalgebra consisting of elements with the part $0$, namely
\begin{equation}
\mathcal W:= \text{Span}_\mathbb C\{a_{-\lambda}\phi_{-\mu}\psi_{-\nu}b_{-\chi} \in \Omega'_1\; |\; -p(\lambda)+p(\mu)-p(\nu)+p(\chi)=0\}.
\end{equation}
Since the Virasoro element $\omega=a_{-1}b_{-1}+\phi_{-1}\psi_{-1}$ is contained in $\mathcal W$,  so $\mathcal W$ is a vertex operator subalgebra.
 Now we will construct a morphism from $\mathcal W$ to $\Omega^{ch}(\mathbb H,\Gamma)/L_1$, sending $v$ to $L(v,1)$, namely
\begin{align*}
L:  \mathcal W &\longrightarrow  \Omega^{ch}(\mathbb H,\Gamma)/L_1\\
v &\longmapsto  L(v,1) 
\end{align*}
where we still use the notation $L(v,1)$ to denote the coset of $L(v,1)$ in $\Omega^{ch}(\mathbb H,\Gamma)/L_1$.
\begin{prop}\label{proposition6.5}
The morphism $L$ defined above is a vertex algebra isomorphism.
\end{prop}
\noindent {\it Proof:} For any $v,w\in \mathcal W$, 
\[
L(v,1)_{(n)} L(w,1)=L(v_{(n)}w,1)+l.o.t
\]
where $l.o.t$ refers to the liftings of nonconstant modular forms which are contained in $L_1$. Hence $L$ is a vertex algebra homomorphism.
According to Lemma \ref{lemma5.1}, 
\[
\{L(a_{-\lambda}\phi_{-\mu}\psi_{-\nu}b_{-\chi},1) \;| \; -p(\lambda)+p(\mu)-p(\nu)+p(\chi)=0\}
\] forms a basis in $\Omega^{ch}(\mathbb H,\Gamma)/L_1$. Thus $L$ is an isomorphism.\qed

In order to understand the maximal quotient $\Omega^{ch}(\mathbb H,\Gamma)/ L_1$, we only need to study $\mathcal W$ due to Proposition \ref{proposition6.5}.
 We first construct a Hermitian form $(\;,\;)$ on $\mathcal W$.  Let $\mathfrak g'$ be a subalgebra of $\mathfrak g$ without $a_0$ and $b_0$, namely
\[
\mathfrak g':= \text{Span} \{ a_{m},b_{m}, \phi_n,\psi_n| m\in \mathbb Z_{\neq 0} ,n\in \mathbb Z   \} \oplus \mathbb C C \subset \mathfrak g,
\]
where $\mathfrak g$ is defined as in (\ref{4.1}).
Let $\alpha$ be an antilinear anti-involution on $\mathfrak g'$ defined by
\begin{gather}
\nonumber \alpha (a_{n})=nb_{-n},\;\;\; \alpha(b_{n})=-\frac{1}{n}a_{-n},  \text{   for  }  n \in \mathbb Z_{\neq 0},\\
\alpha(\psi_{n})= \phi_{-n},\;\;\; \alpha(\phi_{n})= \psi_{-n},\;\;\; \alpha(C)=C, \text{  for  } n\in \mathbb Z.
\end{gather}

Note that there is a natural action of $\mathfrak g'$ on $\Omega'_1$ with the center $C$ acting as the identity operator. $\Omega'_1$ admits a unique positive definite Hermitian form $(\;,\;)$ with the property that $(1,1)=1$,
and for any $x\in \mathfrak g', u,v\in \Omega'_1$
\[
( x u,v  )=(u, \alpha(x) v).
\]
Then the collection 
\[\{a_{-\lambda}\phi_{-\mu}\psi_{-\nu} b_{-\chi}  \in \Omega'_1 | \text{ for all four-tuples } (\lambda,\mu,\nu,\chi)    \}
\] forms an orthogonal basis with respect to $(\;,\;)$. And we obtain a positive definite Hermitian form $(\;,\;)|_{\mathcal W} $ on $\mathcal W$ by restriction, which will be still denoted by $(\;,\;)$. 
We apply the approach in \cite{L} to show the simplicity of $\mathcal W$.
 Let $I$ be a proper ideal  in $\mathcal W$.  Using the Vandemonde matrix, it is not hard to show that $I$ is homogeneous, and hence $I=\sum_{n=0}^\infty I_{n}$, where $I_{n}$ is spanned by the conformal weight $n$ elements in $I$. Obviously $I$ cannot contain the vacuum vector $1$, otherwise $I$ would be the full space $\mathcal W$.
And since the conformal weight $0$ part of $\mathcal W $ is $\mathbb C 1 $,  so $I_{0}=0$, and thus $(1,I)=0$. Therefore for any $w\in \mathcal W$,
\[
(w,I)=(1,\alpha(w) I)=0,
\]
which implies $I$ must be trivial and hence $\mathcal W$ is simple.

\begin{theorem} 
  The vertex operator algebra $\Omega^{ch}(\mathbb H,\Gamma)/L_1$ is simple.
\end{theorem}

\
\

\section{Hecke Operators}
In this section, we will introduce the Hecke operators on the vertex algebra $\Omega^{ch}(\mathbb H,\Gamma)$.
We first extend the action formula in (\ref{3.3}) to the group $GL(2,\mathbb R)_{>0}$. We denote by
\begin{equation} \label{7.1}
g=\begin{pmatrix} \alpha & \beta \\ \gamma  & \delta\end{pmatrix}  \in GL(2,\mathbb R)_{>0},
\end{equation}
and define $\pi(g)=\pi(g')$, where 
 \[
 g'=(\det g)^{-1/2} g \in SL(2,\mathbb R).
 \]
 Then the action formula on the generators can be written as follows.
\begin{align} 
\nonumber \pi(g) a&=a_{-1}  \dfrac{(\gamma b+\delta)^2}{\alpha\delta-\beta\gamma} +2 \dfrac{ \gamma (\gamma b+\delta)}{\alpha\delta-\beta \gamma} \phi _0 \psi_{-1},\\
\nonumber\pi(g) b_{-1}&=b_{-1} \dfrac{\alpha\delta-\beta\gamma}{ (\gamma b+\delta)^{2}},\\
\label{7.2}\pi(g) \psi &=\psi_{-1} \dfrac{(\gamma b+\delta)^{2}}{\alpha\delta-\beta\gamma},\\
\nonumber\pi(g) \phi &=\phi_0 \dfrac{\alpha\delta-\beta\gamma}{(\gamma b+\delta)^{2}},\\
\nonumber\pi(g)f(b)&=f(gb)=f(\dfrac{\alpha b+\beta}{\gamma b +\delta}).
\end{align}

For any integer $k$ and $g$ as in (\ref{7.1}),
we define the weight $k$ operator $[g]_k$ on the ring of holomorphic functions $\mathcal O(\mathbb H)$ by
\[
(f[g]_k)(\tau) =(\det g)^{k/2} (\gamma \tau +\delta)^{-k} f(g\tau),
\]
which gives a right $GL(2,\mathbb R)_{>0}$-action on $\mathcal O(\mathbb H)$.

Let $g\in GL(2,\mathbb Q)_{>0}$, and
\begin{equation} \label{7.3}
\Gamma g \Gamma=\bigcup _{i=1}^n \Gamma g_i,\;\;\; g_i\in GL(2,\mathbb Q)_{>0},
\end{equation}
be a decomposition of $\Gamma g \Gamma$ into right cosets. We recall the double coset operator $[\Gamma g\Gamma]_k$ on $M_k(\Gamma)$ defined by (cf. \cite{Kn} p.273-282)
\[
f [\Gamma g \Gamma]_k :=\sum_{i=1}^n f [g_i]_k.
\]
Similarly we define the double coset operator $\pi(\Gamma g \Gamma)$ on $\Omega^{ch}(\mathbb H,\Gamma)$ by 
\begin{equation}\label{7.4}
\pi(\Gamma g\Gamma) v :=\sum_{i=1}^n \pi(g_i) v, \;\;\; \text{ for any } v\in \Omega^{ch}(\mathbb H, \Gamma).
\end{equation}
Obviously the operator $\pi(\Gamma g \Gamma)$ is independent of the choice of representatives $g_i$, since $v$ is invariant under $\Gamma$-action.  And $\pi(\Gamma g \Gamma) v$ is invariant under the action of $\Gamma$, indeed for $h\in \Gamma$, $g_ih=h_i g_{\sigma(i)}$, where $h_i\in \Gamma$ and $\sigma$ is a permutation of $\{1,2,\cdots,n\}$. So we have
\begin{align*}
\pi(h)\sum_{i=1}^n \pi(g_i)v&=\sum_{i=1} ^n \pi(g_i h) v\\
&=\sum _{i=1}^n \pi(h_i g_{\sigma(i)})v\\
&= \sum_{i=1}^n \pi(g_i)v
\end{align*}
where the last equality is because $\pi$ is a right action and $v$ is invariant under $\Gamma$-action. Hence the operator $\pi(\Gamma g \Gamma)$ is well defined on $\Omega^{ch}(\mathbb H,\Gamma)$.

Denote by $M(n)$ the set of $2\times 2$ matrices over $\mathbb Z$ with determinant $n$. And $M(n,N)$ is a subset of $M(n)$ defined by 
\[
M(n,N):= \left\{\begin{pmatrix} \alpha & \beta \\ \gamma & \delta\end{pmatrix} \in M(n)\; | \; \gamma \equiv 0 \mod N, \text{ gcd}(\alpha,N)=1\right\}.
\]
Now we take our congruence subgroup to be 
\begin{equation}
\Gamma_0(N) :=\left\{ \begin{pmatrix} \alpha & \beta \\ \gamma  & \delta\end{pmatrix} \in SL(2,\mathbb Z)\; | \;  \gamma\equiv 0 \mod N  \right\}.
 \end{equation}
 Then $M(n,N)$ is closed under right and left multiplication by elements in $\Gamma_0(N)$,
and it has a right coset decomposition 
\begin{equation} \label{7.6}
M(n,N)=\bigcup_{i=1}^m \Gamma_0(N) \alpha_i,\;\;\;  \alpha_i \in M(n).
\end{equation}
The Hecke operator $T_k(n)$ on $M_k(\Gamma_0(N))$ and $ \Omega^{ch}(\mathbb H, \Gamma_0(N))$ is defined as follows (cf. \cite{Kn} p.273-282)
\begin{align}\label{7.7}
T_k(n)f:&=n^{\frac{k}{2}-1}\sum_{i=1}^m f [\alpha_i]_k=n^{\frac{k}{2}-1} \sum f[\Gamma_0(N) g \Gamma_0(N)  ]_k,\;\;\;  \text{ for } f\in M_k(\Gamma_0(N)),\\ \label{7.8}
T_k(n) v:&= n^{\frac{k}{2}-1} \sum_{i=1}^m \pi(\alpha_i) v= n^{\frac{k}{2}-1} \sum  \pi(\Gamma_0(N) g\Gamma_0(N)) v,\;  \text{ for } v\in \Omega^{ch}(\mathbb H, \Gamma_0(N))
\end{align}
where the second summations in (\ref{7.7}) and (\ref{7.8}) are over all double cosets of $\Gamma_0(N)$ in $M(n,N).$
Hence the Hecke operator $T_k(n)$ is well-defined.

\begin{prop}\label{prop7.1}
For $g\in GL(2,\mathbb Q)_{>0}$, and $f\in M_{2k}(\Gamma)$, 
\begin{equation} \label{7.9}
\pi(\Gamma g \Gamma) L(w,f)=L(w, f[\Gamma g \Gamma]_{2k}).
\end{equation}
Moreover when $\Gamma=\Gamma_0(N)$, we have
\begin{equation}
T_{2k}(n) L(w,f) =L(w, T_{2k}(n) f).
\end{equation}
\end{prop}

\noindent {\it Proof:} Assume that
\begin{equation} \label{7.11}
\Gamma g\Gamma =\sum_{i=1}^n \Gamma g_i.
\end{equation} 
Denote by $g_i=\begin{pmatrix} \alpha_i & \beta_i\\ \gamma_i&\delta_i \end{pmatrix}= \det g_i ^{1/2}g_i' $, where $g_i'  \in SL(2,\mathbb Q)$.
First, we will prove (\ref{7.9}) for any nonconstant modular form $f \in M_{2k}(\Gamma)$ with $k\geq 1$. According to the proof in Section \ref{Section 4.1}, we may write the lifting as $Af$.
\begin{align*}
\pi(g_i) Af= \pi(g_i') Af&=\pi(g_i') A \pi(g_i')^{-1} \pi(g_i') f\\
&= A  (\det g_i)^{k}  (\gamma_i b+\delta_i)^{-2k} f(g_ib)\\
&= A(  f [ g_i]_{2k}).
\end{align*}
Take summation for all $i$, and hence (\ref{7.9}) holds when $f$ is a nonconstant modular form.

Then we will prove the case for the lifting of the constant modular form. Assume that $\Gamma=\Gamma(1)$ in (\ref{7.11}), and $a_{-\lambda} \phi_{-\mu} \psi_{-\nu} b_{-\chi} +AE$ is a lifting of $1$ in $\Omega^{ch}(\mathbb H,\Gamma(1))$ as in Section \ref{Section 4.2}. Then we have
\begin{align*}
\pi(g_i) (a_{-\lambda} \phi_{-\mu} \psi_{-\nu} b_{-\chi} +AE) =&\pi(g_i')(a_{-\lambda} \phi_{-\mu} \psi_{-\nu} b_{-\chi} +AE)\\
=& a_{-\lambda} \phi_{-\mu} \psi_{-\nu} b_{-\chi}  -A \gamma_i (\gamma_i b+\delta_i)^{-1} +A (\det g_i)  (\gamma_i b+\delta_i)^{-2}   E(g_i b)\\
=&a_{-\lambda} \phi_{-\mu} \psi_{-\nu} b_{-\chi} + A(E[g_i]_2-\gamma_i(\gamma_i b+\delta_i)^{-1})
\end{align*}
Taking summation for all $i$ of the above formula, we have
\begin{align*}
T(\Gamma(1) g\Gamma(1))& (a_{-\lambda} \phi_{-\mu} \psi_{-\nu} b_{-\chi} +AE)\\
&=n a_{-\lambda} \phi_{-\mu} \psi_{-\nu} b_{-\chi}+A \left(\sum_{i=1}^n(E[g_i]_2-\gamma_i(\gamma_i b+\delta_i)^{-1})\right)\in \Omega^{ch}(\mathbb H,\Gamma(1)).
\end{align*}
We subtract the above formula by the lifting of the constant $n$, namely $na_{-\lambda} \phi_{-\mu} \psi_{-\nu} b_{-\chi} +nAE$, then the resulting  element
\[
A \left(\sum_{i=1}^n(E[g_i]_2-\gamma_i(\gamma_i b+\delta_i)^{-1})  -nE\right) 
\]
is still contained in $\Omega^{ch}(\mathbb H,\Gamma(1))$. Denote by $Y$ the holomorphic function 
\[
\sum_{i=1}^n(E[g_i]_2-\gamma_i(\gamma_i b+\delta_i)^{-1})  -nE .
\] 
Since $AY$ is fixed under the action of $\Gamma(1)$,  we have
\[
(\gamma b+\delta)^{-2} Y(hb)=Y(b),\;\;\;\text{ for any } h=\begin{pmatrix} \alpha & \beta \\ \gamma & \delta \end{pmatrix}  \in \Gamma(1),
\]
which implies $Y\in M_{2}(\Gamma(1))=0$. Hence 
\[
T(\Gamma(1) g\Gamma(1))(a_{-\lambda} \phi_{-\mu} \psi_{-\nu} b_{-\chi} +AE)=n(a_{-\lambda} \phi_{-\mu} \psi_{-\nu} b_{-\chi} +AE),
\]
which equals $L(a_{-\lambda} \phi_{-\mu} \psi_{-\nu} b_{-\chi}, 1[\Gamma(1) g \Gamma(1)]_0)$.
Therefore (\ref{7.9}) holds.\qed

We define a group action $\pi'$ of $GL(2,\mathbb Q)_{>0}$ on $\mathcal O(\mathbb H)$ from the right by 
\[
\pi'(g) f(\tau)=( f[g]_2) (\tau) +\frac{6i}{\pi} \gamma(\gamma b+\delta), \;\;\; \text{ for } g=\begin{pmatrix} \alpha & \beta \\ \gamma & \delta \end{pmatrix}\in GL(2,\mathbb Q)_{>0}.
\]
Let $\{g_1,\cdots,g_n\}$ be a complete right coset representatives of $\Gamma(1)$ in $\Gamma(1) g\Gamma(1)$, then we define an operator $T'_g$ on $\mathcal O(\mathbb H)$ by
\[
T'_g f(\tau) :=\frac{1}{n}\sum_{i=1}^n\pi'(g_i) f(\tau).
\]
 From the proof of Proposition {\ref{prop7.1}}, we have
\begin{cor} For any $g\in GL(2,\mathbb Q)_{>0}$, the Eisenstein series $E_2(\tau)$ is a fixed point under $T_g'$-action,
\[
T'_g E_2(\tau)  =E_2 (\tau).
\]
\end{cor}

\section{Acknowledgements}
The paper is based on the auther's dissertation submitted to the Hong Kong University of Science and Technology in 2020. 
The auther wishes to thank his advisor Prof. Yongchang Zhu, for advice. And the auther also wishes to thank Prof. Bailin Song, for discussion.

\end{document}